\documentclass{article}

\usepackage{fullpage}
\usepackage{amsmath,amssymb,amsthm}
\usepackage{bbm} 
\usepackage{bm} 

\usepackage{url}

\usepackage{hyperref}

\usepackage[dvipsnames]{xcolor}

\usepackage{tabularx}                   

\newcommand{\bb}[1]{\mathbb{#1}}

\newcommand{\mc}[1]{\mathcal{#1}}
\newcommand{\msf}[1]{\mathsf{#1}}
\newcommand{\mbf}[1]{\mathbf{#1}}

\newcommand{\fun}{\text{Fun}}
\newcommand{\link}{\msf{Link}}
\newcommand{\bl}{\msf{Bl}}

\newcommand{\gr}{\msf{Gr}}
\newcommand{\spn}{\msf{span}}
\newcommand{\pr}{\msf{pr}}
\newcommand{\lhom}{\text{Hom}^\text{lin}}
\newcommand{\graph}{\msf{Graph}}
\newcommand{\nchoosek}{\textstyle\genfrac{\lbrace}{\rbrace}{0pt}{}}
\newcommand{\rot}{\msf{Rot}}
\newcommand{\trot}{\widetilde{\msf{Rot}}}
\newcommand{\R}{\msf{R}}
\newcommand{\uno}{\mathbbm{1}}
\newcommand{\card}{\msf{card}}
\renewcommand{\dim}{\msf{dim}}
\newcommand{\redc}{C^\text{red}}
\newcommand{\mat}{\msf{Mat}}
\newcommand{\col}{\msf{col}}
\newcommand{\cone}{\msf{Cone}}
\newcommand{\inc}{\msf{inc}}
\newcommand{\ul}[1]{\underline{#1}}
\newcommand{\In}{\msf{In}}
\newcommand{\Out}{\msf{Out}}
\newcommand{\adm}{\msf{adm}}
\newcommand{\Sub}{\msf{Sub}}

\newcommand{\coker}{\msf{coker}}
\renewcommand{\ker}{\msf{ker}}
\newcommand{\op}{\msf{op}}

\def\cP{\mathcal P}

\def\RR{\mathbb R}

\def\ZZ{\mathbb Z}

\def\sC{\mathsf C}
\def\sH{\mathsf H}

\def\sR{\mathsf R}

\newcommand{\lag}{\langle}
\newcommand{\rag}{\rangle}

\newcommand{\ov}{\overline}

\newcommand{\xra}{\xrightarrow}

\newcommand{\e}{\mbf{e}}
\newcommand{\Z}{\mathbb{Z}}

\usepackage{mathtools}                          

\usepackage{float}                                  

\usepackage[normalem]{ulem}                                   

\usepackage{graphics,graphicx,tikz,tikz-cd}                      

\newcommand{\todo}[1]{{\color{red}\textbf{TODO: }{#1}}}

\tikzstyle{vertex}=[circle,fill=black!25,minimum size=20pt,inner sep=0pt]

\newcommand*{\TakeFourierOrnament}[1]{{%
\fontencoding{U}\fontfamily{futs}\selectfont\char#1}}
\newcommand*{\danger}{\TakeFourierOrnament{66}}

\newtheorem{theorem}{Theorem}[section]
\newtheorem{proposition}[theorem]{Proposition}
\newtheorem{corollary}[theorem]{Corollary}
\newtheorem{lemma}[theorem]{Lemma}

\theoremstyle{definition}
\newtheorem{definition}[theorem]{Definition}
\newtheorem{example}[theorem]{Example}

\newtheorem{remark}[theorem]{Remark}

\newtheorem{observation}[theorem]{Observation}

\newtheorem{notation}[theorem]{Notation}

\title{The cohomology of real Grassmannians via Schubert stratifications}
\author{Eric Berry and Scotty Tilton}

\begin{document}
\maketitle

\begin{abstract}

    In this paper, we present a closed formula for the cohomology of real Grassmannians.
    To achieve this, we use a theory of stratified spaces to compute the differentials in a chain complex that computes the cohomology.
    Specifically, we organize Schubert cells as a conically smooth stratified space in the sense of Ayala, Francis, Tanaka — the links therein yield the sought differentials, using methods in differential topology.  
    Further, we identify the isomorphism type of this chain complex and we use this result to
    provide a closed formula for the additive structure of the 
    cohomology of real Grassmannians with arbitrary coefficients.
\end{abstract}

\section{Introduction}

As the moduli space of $k$-dimensional subvector spaces of $\bb{R}^n$, the real Grassmannians, 
$\gr_k(\bb{R}^n)$, play an important role in the study of manifolds. Namely, we can interpret
$H^*(\gr_k(\bb{R}^n))$ as measuring obstructions to natural geometric questions about manifolds.
For instance, let $M \subset \bb{R}^n$
be a $k$-dimensional submanifold of $\bb{R}^n$. Consider the Gauss map
$$\tau_M: M \to \gr_k(\bb{R}^n)~, \qquad x \mapsto T_xM~,$$
that sends a point $x \in M$ to the tangent space of $M$ at $x$. This induces a map on the level of cohomology
$$\tau_M^*:H^*(\gr_k(\bb{R}^n)) \to H^*(M)~.$$ 
For $M$ a compact, orientable submanifold, there exists an element 
$e \in H^k(\gr_k(\bb{R}^n))$ such that $\tau_M(e) = 0$ if and only if $M$ admits a non-vanishing vector field
\cite{characteristic}.\\

To aid in the computation of the cohomology of $\gr_k(\bb{R}^n)$, there is a natural CW structure that one can place on the Grassmannian called the Schubert CW structure \cite{characteristic}. The complex Grassmannian, $\gr_k(\bb{C}^n)$, can also
be endowed with this structure. Since each cell in the Schubert CW strucutre of $\gr_k(\bb{C}^n)$ is even
dimensional, the boundary maps in the corresponding chain complex are zero. Thus, $H^i(\gr_k(\bb{C}^n))$
is freely generated by the Schubert cells of dimension $i$ \cite{griffiths}. A similar situation occurs in the 
cohomology of $\gr_k(\bb{R}^n)$ with $\bb{Z}/2\bb{Z}$ coefficients. Namely, all the boundary maps in the
chain complex are zero. So again, $H^i(\gr_k(\bb{R}^n);\bb{Z}/2\bb{Z})$ is freely generated by the
Schubert cells of dimension $i$ \cite{characteristic}.
Borel computed the rational cohomology algebra of the odd dimensional real Grassmannians \cite{borel}, and
Takeuchi computed the rational cohomology algebra of the even dimensional real Grassmannians \cite{takeuchi}.
More recently, equivariant versions of the rational cohomology have been considered in \cite{carlson}, \cite{he},
\cite{sadykov}. \\

We provide a complete description of the additive structure of the integral cohomology of $\gr_k(\bb{R}^n)$.
Fix $0 \leq k \leq n$, and consider the poset 
$\nchoosek{n}{k} := \fun^{\textsf{inj}} \bigl( \{1< \cdots < k\}, \{1< \cdots < n\}\bigr)$ of injective functors.
For each $0\leq r \leq k$, consider the map
$$d_r : \nchoosek{n}{k} \xrightarrow{\{s_1<\cdots<s_k\} \mapsto\sum_{r\leq i \leq k} s_i-i} \ZZ_{\geq 0}~.$$
In Section \ref{sec:schubert}, we prove the following:
\begin{lemma}\label{introLem1}
    Let $R$ be a commutative ring. The $R$-valued cohomology $\sH^\bullet\bigl( \gr_k(n) ; R\bigr)$ is
    isomorphic with the cohomology of the chain complex $\Bigl(R \left\langle \nchoosek{n}{k} \right\rangle, 
    \delta\Bigr)$ over $R$ whose underlying graded $R$-module is free on the graded set 
    $\nchoosek{n}{k} \xrightarrow{d_1} \bb{Z}_{\geq 0}$ and whose differential evaluates as
    \begin{equation}
        \delta : S = \{s_1< \cdots < s_k\} ~\mapsto~ \sum\limits_{r \in \left\{1 \leq r \leq k ~\mid~ s_{r+1}-s_r >1 
            \text{ and } k-s_r \text{ is odd}\right\}} (-1)^{d_{r-1}(S)} 2 \cdot S_r~,
    \end{equation}
    where $S_r:= \{s_1 < \cdots < s_{r-1} < s_r +1 < s_{r+1} <\cdots< s_k \} \in \nchoosek{n}{k}$.
\end{lemma}

In Section \ref{sec:cohomology} we use Lemma \ref{introLem1} to provide the following remarkable 
decomposition of $C_*(\gr_k(n);\bb{Z})$:
\begin{theorem}\label{introThm}
    Let $R$ be a commutative ring. There is an isomorphism of chain complexes
    \begin{align*}
        C_*(\gr_k(n);\ZZ) &\cong \bigoplus_{S \in \nchoosek{n}{k}_\Out}
            \cone(\ZZ \xrightarrow{2} \ZZ)^{\otimes \In(S)}[d(S)-\card(\In(S))] \\
        &\cong \left(\bigoplus_{S \in \nchoosek{n}{k}_\Out \setminus \nchoosek{n}{k}_\In} \ZZ[d(S)]\right) \oplus
            \left(\bigoplus_{S \in \nchoosek{n}{k}_\msf{Min}} \cone(\ZZ \xrightarrow{2}\ZZ)[d(S)-1]\right)
    \end{align*}
    where $\nchoosek{n}{k}_\Out,\nchoosek{n}{k}_\In, \nchoosek{n}{k}_\msf{Min} \subset \nchoosek{n}{k}$, and $\In(S) \subset \{1, \dots, k\}$.
\end{theorem}

A consequence of Theorem \ref{introThm} is the following closed formula for the 
$R$-cohomology of $\gr_k(n)$:
\begin{corollary}\label{introCor1}
    There is an isomorphism of graded $R$-modules
    $$H^*(\gr_k(n);R) \cong \bigoplus_{S \in \nchoosek{n}{k}} V_S[d_1(S)]~,$$
    where 
    $$V_S := \begin{cases}
                            R, &\text{ if } \msf{In}(S)  = \emptyset = \msf{Out}(S) \\
                            \ker(R\xrightarrow{2}R), &\text{ if } \msf{Min}(\In(S) \cup \Out(S)) \in \Out(S) \\
                            \coker(R\xrightarrow{2}R), &\text{ if } \msf{Min}(\In(S) \cup \Out(S)) \in \In(S)
                    \end{cases},$$
    for $\In(S), \Out(S) \subset \{1, \dots, k\}$.
\end{corollary}

An immedaite consequence of Corollary \ref{introCor1} is the well known result:
\begin{corollary}
    Let $R$ be a commutative ring such that $R\xrightarrow{2}R$ is the zero map, for instance $R=\ZZ/2\ZZ$,
    then 
    $$H^*(\gr_k(n);R) \cong \bigoplus_{S \in \nchoosek{n}{k}} R[d_1(S)]~.$$
\end{corollary}

The difficulty in computing the integral cohomology of 
real Grassmannians lies in the computation of 
the attaching maps in the Schubert CW structure. 
Using $\bb{Z}/2\bb{Z}$ coefficients, one is able to bypass this difficulty,
but to understand the integral cohomology, we must have a clear understanding of the attaching maps. 
In \cite{ck13}, the authors apply a
combinatorial approach via Young diagrams, and
in \cite{jungkind}, Jungkind provides a closed formula for the differentials in the CW chain complex of 
$\gr_k(n)$ with integer coefficients. Using our framework, we recover this formula as Lemma
\ref{introLem1} above. \\

We now outline our approach, a novelty of which
is that our general framework is broadly applicable to a large class of spaces.
We approach this problem of understanding the attaching maps
from the perspective of the recently developed theory of stratified spaces of Ayala-Francis-Tanaka \cite{aft1}. 
Stratified topological spaces are topological spaces equipped with some extra organizational structure called a stratification. By imposing some extra regularity, one arrives at the notion of a conically
smooth stratified space. Conically smooth stratified spaces encompass Euclidean spaces and manifolds, but also
capture manifolds with corners and singularities. 
The extra regularity of a conically smooth stratified space allows us to observe
that the attaching maps are better than just continuous maps to a lower skeleton of a CW structure; these attaching
maps are constructible bundles. In particular, differential topology can be employed for computing the degrees of
maps, which amounts to computing determinants of matrices.
In this paper, we show how the additional structure of a 
conically smooth stratification allows one to compute the homology of a space. We then apply this to the case of
real Grassmannians, $\gr_k(n)$, and provide a description of the chain complex. Further, we explicitly 
determine the additive structure of the cohomology of this complex, and we provide some computations.\\



This paper is organized as follows. 
In Section \ref{sec:grassmannians}, we define the real Grassmannian $\gr_k(n)$ as a smooth manifold. Additionally,
we describe the natural Schubert CW structure on $\gr_k(n)$, and set up the problem of computing the boundary
maps of the corresponding CW chain complex. In Section \ref{sec:stratifications}, we define
stratified topological spaces and conically smooth stratified spaces, and provide a number of examples. 
Further, we prove that one can use a conically smooth structure to compute the homology of a space. 
Section \ref{sec:schubert} introduces the Schubert stratification of $\gr_k(n)$, and then uses the theory of
Section \ref{sec:stratifications} to explicitly compute the boundary maps in the Schubert CW chain complex of
$\gr_k(n)$. Using this formula for the boundary maps, in Section \ref{sec:cohomology}, we establish a
closed formula for the $R$-cohomology of $\gr_k(n)$. Finally, in Section \ref{sec:computations}, we provide
some explicit computations.

\subsection{Acknowledgements}
Both authors would like to acknowledge David Ayala for his constant help and motivation for the project. Without his guidance and patience, we would not have been able to complete this. We would also like to thank Mark Poston who helped in the learning process. ST would like to give a huge thanks to the Undergraduate Scholars Program at Montana State University for funding this research. 
EB acknowledges financial support from NSF-DMS award numbers 1812055 and 1945639, in addition to financial support from the European Research Council (ERC) under the European Union's Horizon 2020 research and innovation program (Grant Agreement No. 768679), as this paper was partially written while EB was visiting IMAG.

\subsection{Notation}
We will now fix some notational conventions that we will use throughout the paper:
\begin{itemize}
    \item For $S$ a finite set, we let $\card(S)$ denote the cardinality of $S$.
    \item For $r\in \bb{Z}_{>0}$, we let $\underline{r}$ denote the set $\{1, \dots, r\}$
    \item For $S$ a set, let $\Sub(S)$ denote the poset of subsets of $S$.
    \item We let $\nchoosek{n}{k}$ denote the set of all cardinality $k$ subsets of the set $\{1, \dots, n\}$.  
    \item Let $\mat_{n \times k}$ denote the collection of $n$-by-$k$ matrices with real coefficients.
    \item Let $\uno_{n\times n}\in \mat_{n \times n}$ denote the $n$-by-$n$ identity matrix.
    \item For each $1\le i \le n$, $\e_i \in \bb{R}^n$ is the $i$-th standard basis vector of $\bb{R}^n$. 
        That is, it is the vector consisting of a 1 in the $i$-th entry, and 0's in all other entries.
    \item For $S \in \nchoosek{n}{k}$, we define $\bb{R}^S$ to be the span
        $$\bb{R}^S := \spn\{\e_i ~|~ i \in S\}~.$$
    \item For a $X$ a topological space, we let $X^+$ denote the set, $X \amalg \{\infty\}$, where we attach a
        single disjoint point to $X$. This set is equipped with the following topology: $U \subset X^+$ is open
        iff $U$ is an open subset of $X$, or $U = (X \setminus C) \amalg \{\infty\}$, for some closed and compact
        $C \subset X$.
    \item For $\ast \in X$ a pointed topological space, we let $\Sigma X$ denote the reduced suspension,
        $$\Sigma X := (X \times [0,1])/\sim~,$$
        where $\sim$ is the equivalence relation defined by $(x,1) \sim (x',1), (x,0) \sim (x',0),$ and 
        $(\ast, t) \sim (\ast,t')$, for all $x,x' \in X$, and $t,t' \in [0,1]$.
    \item We use homological grading conventions: $A[r]_n := A_{r+n}$.
\end{itemize}

\section{Grassmannians}
\label{sec:grassmannians}
Grassmannians are the moduli space
of subspaces of a fixed vector space. In this paper, we restrict our attention
to subspaces of Euclidean space. These are the \emph{real} Grassmannians:

\begin{definition}
    Let $n,k$ be positive integers such that $k \le n$. We define the 
    \emph{Grassmannian}, $\gr_k(n)$, to be the set consisting of all $k$-dimensional subspaces of $\bb{R}^n$.
\end{definition}

The simplest case to understand is when $n=1$, as $\gr_1(1)$ is just a singleton. 
In general, $\gr_n(n)$ is also a singleton. The first nontrivial example is $(n,k) = (2,1)$.
One way to think about this set is to observe that, by taking its span, each choice of a unit vector
in $\bb{R}^2$ defines a 1-dimensional subspace, hence a point in $\gr_1(2)$.
Further, by observing that antipodal vectors define the same subspace, 
we see that the following map of sets is a bijection
$$[0,\pi) \xrightarrow{\cong} \gr_1(2), \qquad \theta \mapsto \spn\{(\cos\theta, \sin\theta)\}~.$$

Thinking of a Grassmannian as a set does not capture its entire being. 
In the above example of $\gr_1(2) \cong [0,\pi)$, if we sweep through the angles $\theta$, as we approach 
$\theta=\pi$, our subspace is getting closer to the $x$-axis, which is represented by the point $0 \in [0,\pi)$. In 
other words, thinking of $\gr_1(2)$ as the set $[0,\pi)$, ignores the fact that points near $\pi$ are actually near $0$
as well. We will make this intuition precise below, and we will see that $\gr_1(2)$ is actually diffeomorphic to 
$S^1$. \\

More generally, each $p \in \gr_k(n)$ can be representated by a $k$-plane in $\bb{R}^n$, 
and we know how to wiggle planes in Euclidean space.
This suggests two things. First, we know which planes are close to a given plane, so we might expect
$\gr_k(n)$ to be a topological space. Second, we know how to budge a plane into a nearby plane, so we might further
expect $\gr_k(n)$ to possess the structure of a manifold. 
This is indeed the case, as we briefly recall in Section
\ref{mfld}. Even better than possessing a manifold structure, there is a natural way of decomposing $\gr_k(n)$
into pieces called \emph{Schubert cells}. These endow $\gr_k(n)$ with the structure of a CW complex. We give
a description of this Schubert CW structure and the corresponding CW chain complex in Section \ref{cell}. 
We refer the reader to \cite{characteristic} for a more in depth treatment of this material.  \\

\subsection{Manifold structure}\label{mfld}
Let us fix positive integers $n$ and $k$, with $k \le n$. We define a manifold structure on $\gr_k(n)$ by 
realizing it as the quotient of an open subset of Euclidean space called the Stiefel space.

\begin{definition}
    We define the \emph{Stiefel space} $V_k(n)$ to be the collection of all injective, linear maps 
    from $\bb{R}^k$ to $\bb{R}^n$:
    $$V_k(n) := \{\bb{R}^k \xrightarrow{\varphi} \bb{R}^n ~|~ \varphi \text{ is injective and linear}\} \subset
        \mat_{n \times k}~.$$
\end{definition}

A point in $V_k(n)$ is a matrix $A \in \mat_{n \times k}$ whose columns are linearly independent.
Linear independence is an open condition, which means that $V_k(n) \subset (\bb{R}^n)^{\times k}$ is an
open subset of Euclidean space. As such we consider $V_k(n)$ as a topological space via the subspace topology.
Each point $A \in V_k(n)$ 
defines a $k$-dimensional subspace of $\bb{R}^n$ by taking its column space. 
This defines a map
$$\col : V_k(n) \to \gr_k(n), \qquad A \mapsto \col(A)~.$$
Observe that this map is surjective: given $V \in \gr_k(n)$, choose a basis for it and define a matrix whose columns
consist of those basis vectors. This defines a point in $V_k(n)$ whose column space is $V$.
Thus, we endow $\gr_k(n)$ with the quotient topology induced by the map $\col$. 
Further, one can show that $\gr_k(n)$ is a compact topological space, and that it has a natural manifold structure:

\begin{proposition}
    The space $\gr_k(n)$ is a compact, smooth manifold of dimension $k(n-k)$. 
\end{proposition}
{\it Proof.} We defer the reader to \cite{characteristic} for the claims of compactness and Hausdorfness. 
We will define a smooth atlas on $\gr_k(n)$ as follows. Let $\pr_S: \bb{R}^n \to \bb{R}^S$
denote the orthogonal projection $x \mapsto \sum_{i \in S} \langle x,e_i\rangle e_i$.
Given $S \in \nchoosek{n}{k}$, consider the subset
$$U_S := \left\{ V \subset \bb{R}^n ~|~ V \xhookrightarrow{\iota} \bb{R}^n \xrightarrow{\pr_S} \bb{R}^S
    \text{ is an isomorphism}\right\} \subset \gr_k(n)~.$$
Define a map
$$\lhom(\bb{R}^S, \bb{R}^{\{1, \dots, n\} \setminus S}) \to U_S$$
whose value on an injective linear map $F$ is its graph:
$$F \mapsto \graph(F) := \{x+F(x) ~|~ x \in \bb{R}^S\} \subset \bb{R}^n~.$$
In fact, this map is a bijection with inverse given by the following:
given $V \in U_S$, define a linear map $F \in \lhom(\bb{R}^S, \bb{R}^{\{1, \dots, n\} \setminus S})$
via the assignment
$$e_i \mapsto \pr_S^{-1}(e_i) \cap V - e_i~,$$
for all $i \in S$.
Further, one can check that this in fact defines a homeomorphism. 
Thus, for each $S \in \nchoosek{n}{k}$,
we have a homeomorphism
$$\varphi_S: \bb{R}^{k(n-k)} \cong \lhom(\bb{R}^S, \bb{R}^{\{1, \dots, n\} \setminus S}) \xrightarrow{\cong}
    U_S~.$$
We leave it as an exercise to check that the collection $\{\bb{R}^{k(n-k)} \xrightarrow{\varphi_S} U_S\}$ 
defines a smooth atlas for $\gr_k(n)$. (See e.g. 7.8 in \cite{manifolds}.) 
\qed \\

\subsection{CW structure}\label{cell}
There is a natural CW structure on Grassmannians: the \emph{Schubert} CW structure. The set of Schubert cells
is in bijection with the set $\nchoosek{n}{k}$: 
\begin{definition}\label{nchoosek}
    For positive integers $k\le n$, let 
    $$\nchoosek{n}{k}:= \{S \subset \{1,\dots,n\} ~|~ \card(S) = k\}~,$$
    denote the set consisting of all cardinality $k$ subsets of the set $\{1,\dots,n\}$.
    Define a partial order on the elements of $\nchoosek{n}{k}$ by declaring 
    $S = \{s_1<\cdots<s_k\}\leq T = \{t_1<\cdots<t_k\}$ to mean $s_i\leq t_i$ for each $i\in\{1,\dots,k\}$.
\end{definition}

\begin{example}
    The poset $\nchoosek{4}{2}$ can be depicted as
    $$
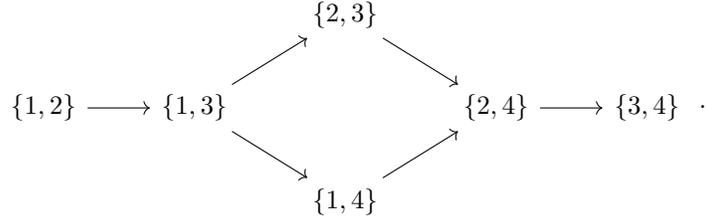
\begin{tikzcd}
            &&\{2,3\}\arrow[dr]&&\\
            \{1,2\}\arrow[r]&\{1,3\}\arrow[ur]\arrow[dr]&&\{2,4\}\arrow[r]&\{3,4\}\\
            &&\{1,4\}\arrow[ur]&&
        \end{tikzcd}~.$$
\end{example}

For $S \in \nchoosek{n}{k}$, the Schubert cell corresponding to $S$
is the subspace
$$\gr_k(n)_S := \{ V \in \gr_k(n) ~|~ 
    \text{$S$ is the maximal element in } {\nchoosek{n}{k}} \text{ for which }V \in U_S \}~.$$
One can check that such a maximum indeed exists, and is unique.
For example, consider $\nchoosek{2}{1} = \{\{1\}<\{2\}\}$, and take $S = \{1\}$. Then $\gr_1(2)_{\{1\}}$
is simply a singleton, because the only 1-dimensional subspace of $\bb{R}^2$ that does not project isomorphically
onto the $y$-axis, $\bb{R}^{\{2\}}=\spn\{e_2\}$, is the $x$-axis, $\bb{R}^{\{1\}} = \spn\{e_1\}$. This also
tells us that $\gr_1(2)_{\{2\}}$ consists of all other points in $\gr_1(2)$. Thus, $\gr_1(2)$ is comprised of a single
0-dimensional cell, 
$$\gr_1(2)_{\{1\}} = \left\{\spn\begin{bmatrix}
    1\\0
\end{bmatrix}\right\}~,$$
and a single 1-dimensional cell, 
$$\gr_1(2)_{\{2\}} = \left\{\spn \begin{bmatrix}
    \ast\\ 1
\end{bmatrix} ~\bigg|~ \ast\in \RR\right\}~,$$
as depicted by the
following picture
$$
    
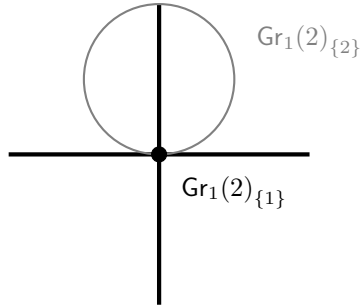
\begin{tikzpicture}
        \draw[ultra thick,black] (-2,0)--(2,0);
        \draw[ultra thick,black] (0,-2)--(0,2);
        \draw[thick, gray] (0,1) circle (1);
        \draw[black, fill = black] (0,0) circle (.1);
        \draw (2,1.5) node {\textcolor{gray}{${\gr_1(2)}_{\{2\}}$}};
        \draw (1,-.5) node {${\gr_1(2)}_{\{1\}}$};
    \end{tikzpicture}~.
$$
To name a CW structure, one must define attaching maps that specify how the higher dimensional cells 
attach to the lower dimensional cells. 
We will now explicitly describe this cellular decomposition of $\gr_k(n)$.\\

For each $S \in \nchoosek{n}{k}$, we define
a map $\rot_S$ from a closed cube to $\gr_k(n)$, whose interior is homeomorphic with $\gr_k(n)_S$, and whose
boundary maps to strictly lower dimensional cells. 
To define $\rot_S$, we will use a product of rotation matrices, which are elements of $O(n)$. Since $O(n)$ is
non-Abelian, the order in which we multiply these matrices is critical. To obtain the correct order of multiplication,
we use a particularly ordered index set, $Z_S$:

\begin{definition}
For $S=\{s_1< \cdots < s_k\} \in \nchoosek{n}{k}$,
        define the set 
        $$Z_S:=\left\{(i,j)\mid i\in\{1,\dots,k\}, j\in \{i,\dots,s_i-1\}\right\}~.$$
        We equip this set with a total ordering by declaring $(i,j) \le (i',j')$ to mean either $i < i'$, or $i=i'$ and
        $j \ge j'$. See Remark \ref{indexSet} for a further explanation of this ordering.
\end{definition}

For $1 \le j < n$ and $\theta \in \bb{R}$, let $\R_j(\theta)$ denote the $n \times n$ matrix implementing
a rotation of $\bb{R}^n$ in the oriented $(j<j+1)$-plane by an angle of $\theta$:
$$\R_j(\theta):=\begin{bmatrix}
    &&&&&&&&\\
    &\uno_{j-1\times j-1} &&&&&&&\\
    &&&&&&&&\\
    &&&\cos(\theta)&-\sin(\theta)&&&&\\
    &&&\sin(\theta)&\cos(\theta)&&&&\\
    &&&&&&&&\\
    &&&&&&&\uno_{(n-j-1)\times(n-j-1)}&\\
    &&&&&&&&  
    \end{bmatrix}~,$$
where the middle block lies in the $j$ and $j+1$ rows and columns, and all unspecified entries are 0.
For example, for $n=3$,
$$\R_2(\theta) = \begin{bmatrix} 1 & 0 & 0 \\ 0 & \cos\theta & -\sin\theta \\ 0 & \sin\theta & \cos\theta 
                                  \end{bmatrix}~.$$
Note that for each $j$, this defines a continuous map 
$$[0,\pi] \to O(n), \qquad \theta \mapsto \R_j(\theta)~.$$
Given a matrix $A \in O(n)$, the first $k$ columns of $A$ is an ordered set of $k$ linearly independent vectors in 
$\bb{R}^n$, which is a point in $V_k(n)$.
Note the map
$$O(n) \to V_k(n), \qquad A \mapsto A\uno_{n\times k}~,$$
where $\uno_{n\times k}$ denotes the first $k$ columns of the $n\times n$ identity matrix. This map is 
evidently continuous, since matrix multiplication is continuous. \\

\begin{remark}\label{indexSet}
    The reason for the above defined ordering on $Z_S$ is to ensure that the matrix multiplication in the 
    definition of $\trot_S$
    occurs in the correct sequence. Namely, we have set this up so that if we evaluate $\trot_S$ at the point
    $(\pi/2)_{Z_S}$, then we obtain the plane $\bb{R}^S$. The defined ordering first moves $\e_k$ into $\e_{s_k}$,
    and ensures that the subsequent rotations will not move 
    $\e_{s_k}$ further. Next, it moves $\e_{k-1}$ into $\e_{s_{k-1}}$, again ensuring that subsequent rotations will
    not move either $\e_{s_k}$ or $\e_{s_{k-1}}$. This process continues, until the final result produces $\bb{R}^S$. 
\end{remark}

These maps, together with the group structure of $O(n)$, and the particularly ordered index set, $Z_S$
(as just explained in Remark \ref{indexSet}), 
allow us to define the composite continuous map
$$\trot_S: \prod_{(i,j) \in Z_S}[0,\pi] \xrightarrow{(\R_j(-))_{(i,j) \in Z_S}} \prod_{Z_S} O(n) \xrightarrow{\text{mult}}
        O(n) \xrightarrow{- \cdot \uno_{n \times k}} V_k(n)~,$$
that sends 
$$(\theta_{(i,j)})_{(i,j) \in Z_S} \mapsto
        \left(\prod_{(i,j) \in Z_S} \R_j(\theta_{(i,j)})\right) \uno_{n\times k}~.$$
Finally, postcomposing with $\col:V_k(n) \to \gr_k(n)$, we obtain the continuous map
$$\rot_S: \prod_{Z_S}[0,\pi] \xrightarrow{\trot} V_k(n) \xrightarrow{\col} \gr_k(n)~.$$

Having defined the attaching maps, we now show that the collection of sets $\gr_k(n)_S$ do form a CW structure
on $\gr_k(n)$.
First,
we establsih a homeomorphism of $\gr_k(n)_S$ with the interior of the closed disk $\prod_{Z_S} [0,\pi]$. 
In particular, this will tell us that the dimension of the cell $\gr_k(n)_S$ is the cardinality of $Z_S$. Let us
denote this dimension by $d(S)$:
$$d(S) := \card(Z_S) = \dim(\gr_k(n)_S)~.$$
\begin{observation}\label{dimension}
    Recall that for $S=\{s_1< \cdots <s_k\} \in \nchoosek{n}{k}$, 
    $$d(S) := \card(Z_S) = \sum_{i \in \{1, \dots, k\}} \card(\{i, \dots, s_i-1\}) = \sum_{i \in \{1, \dots, k\}} s_i-i~.$$
\end{observation}

\begin{lemma}\label{homeo}
    For $S \in \nchoosek{n}{k}$, the restriction of $\rot_S$ to the interior of $\prod_{Z_S}[0,\pi]$~,
    $$\rot_S : \prod_{Z_S} (0,\pi) \to \gr_k(n)~,$$
    is a homeomorphism onto $\gr_k(n)_S$.
\end{lemma}
{\it Proof.} We prove this by induction on $k$. 
First, assume $k=1$, and let $S=\{s\} \in \nchoosek{n}{1}$. Note that for $V \in \gr_1(n)_S$, 
$$\col^{-1}(V)=\left\{\mbf{x}=(x_1, x_2, \dots, x_s, 0, \dots, 0) ~|~ \sum x_i^2=1, ~x_s \neq 0, ~\col(\mbf{x}) \cong V\right\}~.$$
Observe there is a unique $\mbf{x} \in \col^{-1}(V)$ such that $x_s >0$. Further, we claim there is a unique element 
$\theta^\mbf{x} \in \prod_{Z_S}(0,\pi)$ such that $\trot_S(\theta^\mbf{x}) = \mbf{x}$. As we now explain, this follows from the
injectivity of each $\sR_j:(0,\pi) \to O(n)$, and the order in which the $\sR_j$ matrices are multiplied, i.e., the order 
of $Z_S$. Note that $Z_S = \{(1,s-1) < \cdots < (1,1)\}$. As explained in Remark \ref{indexSet}, this ordering on
$Z_S$ is such that the standard basis vector $\e_1$ gets rotated into $\mbf{x}$ by first rotating in the $(1,2)$-plane
by an angle $\theta_{(1,1)} \in (0,\pi)$, and then rotating in the $(2,3)$-plane by an angle $\theta_{(1,2)} \in (0,\pi)$,
etc., until a final rotation in the $(s-1,s)$-plane by an angle $\theta_{(1,s-1)} \in (0,\pi)$. Since a rotation in
the $(j,j+1)$-plane does not effect the $j-1$-coordinate, we see that for each $j \in \{1, \dots, s-1\}$, there is a 
unique value $\theta^\mbf{x}_{(1,j)} \in (0,\pi)$ such that the first $j$ coordinates of the vector
$$\sR_j(\theta^\mbf{x}_{(1,j)})\sR_{j-1}(\theta^\mbf{x}_{(1,j-1)})\cdots \sR_1(\theta^\mbf{x}_{(1,1)})\e_1$$ 
are equal to the first $j$ coordinates of $\mbf{x}$. This completes the base case.\\

Now, assume $k>1$. For $S=\{s_1< \cdots<s_{k+1}\} \in \nchoosek{n}{k+1}$, let $\widehat{S}$ denote the set
$S \setminus \{s_{k+1}\}$. Note that for $V \in \gr_{k+1}(n)_S$,
$$\col^{-1}(V) = \left\{ \left. A=\begin{bmatrix} a_{1,1} & \cdots & a_{1,\ell} & \cdots & a_{1,k+1} \\
                                                                                        \vdots &  &  & & \\
                                                                                        a_{s_1,1} & & \vdots & &  \\
                                                                                        0 & &  & & \vdots\\
                                                                                         & & a_{s_\ell,\ell} & & \\
                                                                                         & & 0 & &  \\
                                                                                         \vdots & & & & a_{s_{k+1},k+1} \\
                                                                                         & & \vdots & & 0\\
                                                                                         & & & & \vdots \\
                                                                                         0 & 0 & 0 & 0 & 0\\
                                                            \end{bmatrix} ~\right\vert~ \mbf{a}_i \cdot \mbf{a}_j =\delta_i^j,~
                                                            \msf{col}(A) \cong V,~a_{s_i,i} \neq 0 \right\}~,$$
where $\mbf{a}_i$ denotes the $i$th column of $A$, and $\delta_i^j$ is the Kronecker delta. Again, there is a
unique $A_V \in \col^{-1}(V)$ such that $a_{s_i,i} >0$ for all $i \in \{1, \dots, k+1\}$.
From the induction hypothesis, we know the map
$$\rot_{\widehat{S}} : \prod_{Z_{\widehat{S}}} (0,\pi) \to \gr_k(n)_{\widehat{S}}$$
is a homeomorphism. By an analysis similar to that of the base case, we can extend this to a homeomorphism
$$\rot_S : \prod_{Z_S} (0,\pi) \to \gr_{k+1}(n)_S$$
by including the extra necessary rotations to move $\e_{k+1}$ to the last column vector of $A_V$. 
In particular, $Z_S \setminus Z_{\widehat{S}} = \{(k+1, s_{k+1}-1), \dots, (k+1,k+1)\}$, so we define
$$\trot_S := \trot_{\widehat{S}} \circ \msf{R}_{s_{k+1}}(\theta_{(k+1,s_{k+1})}) \cdots 
    \msf{R}_{k+1}(\theta_{(k+1,k+1)})~,$$
for some choice of angles. As before, there is a unique choice of angles in $(0,\pi)$ to achieve these rotations, and
the addition of these rotations will only effect $\e_{k+1}$.
\qed \\

Next, we establish that the boundaries of these cells attach to strictly lower dimensional cells. In what follows, 
denote
$$\partial \prod_{Z_S} [0,\pi] := \prod_{Z_S} [0,\pi] \setminus \prod_{Z_S} (0,\pi)~.$$
\begin{lemma}\label{boundary}
    For $S \in \nchoosek{n}{k}$, the restriction of $\rot_S$ to the boundary factors through the union of the lower
    dimensional cells:
    $$\begin{tikzcd}
        \partial \displaystyle\prod_{Z_S} [0,\pi]  \arrow[r,hook] \arrow[dr,dashed] & \displaystyle\prod_{Z_S}[0,\pi]  
            \arrow[r,"\rot_S"] & \gr_k(n) \\
        & \displaystyle\bigcup_{T <S} \gr_k(n)_T \arrow[ur, hook] & 
        \end{tikzcd}~.$$
\end{lemma}
{\it Proof.} Let $\Theta=(\theta_{(i',j')} \in \partial \displaystyle\prod_{Z_S} [0,\pi]$. 
We will show that for each $1 \le i \le k$, and each $j>s_i$, the $(j,i)$ entry of $\trot_S(\Theta)$ is 0. Fix such a
$(j,i)$ pair. Let 
$$B:= \prod\limits_{(i',j') \in Z_S,~i'\le i} R_{j'}(\theta_{(i',j')})~,$$
and
$$A := \prod\limits_{(i',j') \in Z_S,~i'>i} R_{j'}(\theta_{(i',j')})~.$$
Since $j\ge s_i+1$, and each $(i',j') \in Z_S$ with $i'\le i$ has $j' <s_i$, the $j$-th row of $B$ is $\e_j^T$. Further,
since each $(i',j') \in Z_S$ with $i'>i$ has $j'\ge i'>i$, the $i$-th column of $A$ is $\e_i$. Finally, since
$i\le s_i < j$, $\e_j^T \cdot \e_i = 0$.
\qed \\

\begin{proposition}
    The collection
    $$\left\{(\gr_k(n)_S, \rot_S) ~|~ S \in \nchoosek{n}{k} \right\}$$
    defines a CW structure on $\gr_k(n)$. \qed
\end{proposition}
{\it Proof.} This follows immediately from Lemmas \ref{homeo} and \ref{boundary}. 
\qed \\

\subsection{The Schubert CW chain complex}
To give a general description of the CW chain complex for the Schubert CW structure on $\gr_k(n)$, 
we need to know how many cells there are of each
dimension. It turns out that these can be counted by a partition function. For positive integers $k,n$ with $k < n$, 
and $1 \le m \le k(n-k)$, define 
$$p_k^n(m) := \left\{\{i_1 \le \cdots \le i_\ell\} ~|~ \sum_{1 \le j \le \ell} i_j = m, ~ \ell\le k, ~ i_j \le n-k\right\}~.$$
to be the set of all partitions of $m$ as the sum of at most $k$ positive
integers, each of which is less than or equal to $n-k$. For consistency, let us define $p_k^n(0) := \{0\}$.
\begin{lemma}\label{numOCells}
    There is a bijection between the set of $m$-dimensional cells of the Schubert decomposition of $\gr_k(n)$ 
    and the set $p_k^n(m)$.
\end{lemma}
{\it Proof.} Consider the map of sets
$$p_k^n(m) \xrightarrow{\varphi} \nchoosek{n}{k}~,$$ 
via the assignment
$$\{i_1 \le \cdots \le i_\ell\} \mapsto \{1,2,\dots,k-\ell,k-(\ell-1)+i_1, k-(\ell-2)+i_2, \dots, k+i_\ell\}~.$$
First note that $\{1,2,\dots,k-\ell,k-(\ell-1)+i_1, k-(\ell-2)+i_2, \dots, k+i_\ell\}$ indeed has cardinality $k$, and is thus
an element of $\nchoosek{n}{k}$. By definition, this map is injective. We will now show that the image of $\varphi$
is the collection of sets $S \in \nchoosek{n}{k}$ such that $d(S) = m$.
Let $S=\{s_1< \cdots< s_k\}$ be such that $d(S) = m$. By Observation \ref{dimension},
$$\sum_{j=1}^k s_j-j = m~.$$
Therefore there exists some set $L\subset \{1,\dots,k\}$ where for each $j\in L$, $s_j-j\neq 0$ and 
$$\sum_{j\in L}s_j-j = m~.$$ 
Then $I:=\{s_j-j ~|~ j\in L\}$ is an element of $p^n_k(m)$ for which $\varphi(I)=S$.
Finally, we show that the only sets in the image of $\varphi$ are those $S \in \nchoosek{n}{k}$ for which $d(S)=m$.
Given $S=\{1,2,\dots, k-\ell,k-(\ell-1)+i_1,\dots,k+i_\ell\} \in \nchoosek{n}{k}$ in the image of 
$\varphi$, by Observation \ref{dimension}, the dimension of the corresponding Schubert cell is 
$$d(S) = \sum_{j=1}^k s_j-j = i_1+\cdots +i_\ell = m~.$$
\qed \\

Lemma \ref{numOCells}
tells us that the CW chain complex for the Schubert CW structure on $\gr_k(n)$ is of the form
$$0 \to \bigoplus_{p^n_k(k(n-k))} \bb{Z} \xrightarrow{\partial_{k(n-k)}} \cdots 
    \xrightarrow{\partial_{j+1}} 
    \bigoplus_{p^n_k(j)} \bb{Z} \xrightarrow{\partial_j} \cdots \xrightarrow{\partial_1} 
    \bigoplus_{p^n_k(0)} \bb{Z}
    \xrightarrow{\partial_0} 0~.$$
Therefore, in order to compute the homology of this chain complex, we need to have a description of these
boundary maps. Proposition \ref{homology}, to come, will provide us with a way of computing these
boundary maps.






\section{Homology of stratified spaces}
\label{sec:stratifications}
As we saw in the preceeding section, the remaining task to compute the homology of $\gr_k(n)$ is to 
determine the differentials. In order to do this, we exploit some natural extra regularity on the Schubert
decomposition of $\gr_k(n)$, namely, 
a stratification. For the purposes of computing homology, we need to impose some further regularity on this stratification, that of a \emph{conically smooth structure}. In this section we give a brief introduction to the theory of stratified spaces 
\`a la \cite{aft1}, including the notion of conical smoothness, and then we
discuss how this theory enables us to compute the homology of conically smooth stratified spaces. \\

\subsection{Stratified spaces}
In this section, we introduce the basic definitions of the theory of stratified spaces that we will use  in the
remainder of the paper. 

\begin{definition}
    A \emph{stratified topologcial space} is a triple $(X \xrightarrow{\psi} \mc{P})$ consisting of
    \begin{itemize}
        \item a paracompact, Hausdorff topological space, $X$,
        \item a poset $\mc{P}$, and
        \item a continuous map $\psi:X \to \mc{P}$, where $\mc{P}$ is equipped with the downward closed topology:
            $U \subset \mc{P}$ is closed if for each $p \in U$, if $q<p$, then $q \in U$.
    \end{itemize}
\end{definition}

\begin{definition}
    For $X \xrightarrow{\varphi} \mc{P}$ a stratified topological space, and $p \in \mc{P}$, the \emph{$p$-stratum} 
    of $X$ is the subspace
    $$X_p :=  \varphi^{-1}(p) \subset X~.$$
\end{definition}

\begin{example}
    Let $M$ be a smooth manifold, and let $W \hookrightarrow M$ be a properly embedded submanifold. The map
    $$M \to \{0<1\}, \qquad m \mapsto \begin{cases} 0, &\quad \text{if } m \in W \\
                                                                                       1, &\quad \text{otherwise}
                                                               \end{cases}$$
    exhibits $M$ as a stratified topological space with strata $M_0 = W$ and $M_1 =M \setminus W$.
\end{example}

\begin{example}
    Let $X$ be a topological space equipped with a CW structure. 
    The skeleta of the CW structure give rise to a stratification of $X$.
    The stratifying poset is the nonnegative integers with their natural partial order, $(\ZZ_{\ge0},\le)$.
    For $k \in \ZZ_{\ge0}$, let $X_k$ denote the $k$-skeleton of $X$. That is, $X_k$ is the union of all cells, 
    $X_\alpha$, in the CW structure of $X$ such that $\dim(X_\alpha) \le k$:
    $$X_k := \bigcup_{\dim(X_\alpha)\le k} X_\alpha~.$$
    Notice that each element $x \in X$ lies in $X_k \setminus X_{k-1}$ for some unique $k$.
    Define the map $X \to \ZZ_{\ge0}$ to be such that $x \in X$ maps to the unique $k \in \ZZ_{\ge0}$ such that
    $x \in X_k \setminus X_{k-1}.$
\end{example}

\begin{definition}
    Let $p \in \mc{P}$ be a minimal element of a poset. The \emph{link of $\mc{P}$ along $p$} is defined to be the
    poset
    $$\link_p(\mc{P}) := \mc{P}_{p<}:= \{q \in \mc{P} ~|~ p<q\}~.$$
\end{definition}

\begin{definition}
    Let $p \in \mc{P}$ be a minimal element of a poset. The \emph{blowup of $\mc{P}$ along $p$} is defined to be
    the poset
    $$\bl_p(\mc{P}) := \link_p(\mc{P}) \times \{0<1\} 
        \coprod_{\link_p(\mc{P}) \times \{1\}} \mc{P} \setminus \{p\}~.$$
\end{definition}

\begin{example}\label{planePoset}
    Consider the poset
    $$\cP = \begin{tikzcd}
            NW&N\arrow[r]\arrow[l]&NE\\
            W\arrow[u]\arrow[d]&O\arrow[u]\arrow[d]\arrow[r]\arrow[l]&E\arrow[u]\arrow[d]\\
            SW&S\arrow[r]\arrow[l]&SE
        \end{tikzcd}~.$$
    The link and blowup of $\cP$ along $O$ are as follows:
    $$\link_O(\cP) = \begin{tikzcd}
            NW&N\arrow[r]\arrow[l]&NE\\
            W\arrow[u]\arrow[d]&&E\arrow[u]\arrow[d]\\
            SW&S\arrow[r]\arrow[l]&SE
        \end{tikzcd}~,\quad
            \bl_O(\cP) = \begin{tikzcd}[row sep= 0.7em,column sep = 0.7em]
            (NW,1)&&(N,1)\arrow[rr]\arrow[ll]&&(NE,1)\\
            &(NW,0)\arrow[ul]&(N,0)\arrow[r]\arrow[l]\arrow[u]&(NE,0)\arrow[ur]\\
            (W,1)\arrow[uu]\arrow[dd]&(W,0)\arrow[l]\arrow[u]\arrow[d]&&(E,0)\arrow[r]\arrow[u]\arrow[d]&(E,1)\arrow[uu]\arrow[dd]\\
            &(SW,0)\arrow[dl]&(S,0)\arrow[r]\arrow[l]\arrow[d]&(SE,0)\arrow[dr]\\
            (SW,1)&&(S,1)\arrow[rr]\arrow[ll]&&(SE,1)\\
        \end{tikzcd}~.$$
\end{example}

\begin{definition}
    Let $X \to \mc{P}$ be a smooth manifold equipped with the additional structure of a stratified topological
    space. Further, assume that each stratum, $X_p$, is a smooth submanifold of $X$. For $p \in \mc{P}$ a minimal element of the stratifying poset, we define the \emph{link of $X$ along $X_p$} to be the 
    stratified topological space with underlying space the unit sphere bundle of the normal bundle
    $$\link_{X_p}(X) := S(N_{X_p \subset X})~,$$
    and stratifying poset given by $\link_p(\mc{P})$. 
\end{definition}

Note that this definition in fact exhibits the link as a smooth fiber bundle
$$\pi:\link_{X_p}(X) \to X_p~.$$


A choice of a tubular neighborhood of $X_p \subset X$ gives us a smooth map 
$N_{X_p \subset X} \hookrightarrow X$. 
We then get a smooth map from the thickened link to $X$:
$$\gamma: \link_{X_p}(X) \times (0,\infty) \hookrightarrow N_{X_p \subset X} \hookrightarrow X~.$$

\begin{definition}
    Let $X \to \mc{P}$ be a smooth manifold equipped with the additional structure of a stratified topological space, and let $p \in \mc{P}$ be a minimal element.
    We define the \emph{blowup of $X$ along $X_p$} to be the pushout
    $$\bl_{X_p}(X):=\left(\link_{X_p}(X) \times [0,\infty)\right) \coprod_{\link_{X_p}(X) \times (0,\infty)} 
    X \setminus X_p~.$$
    This space is naturally stratified by the poset $\bl_p(\mc{P})$.
\end{definition}

\begin{example} 
    We can stratify $\RR^2$ by the poset $\cP$ from Example \ref{planePoset}. We specify the structure map
    $\RR^2 \to \cP$ in Figure \ref{fig:plane} by labeling the strata. There is a single 0-dimensional
    stratum $\RR^2_O$ consisting of the origin.
    The resulting blowup and link of $(\RR^2\to \cP)$ along $\RR^2_O$ are given in Figure~\ref{fig:plane}.
\end{example}
\begin{figure}
    \begin{center}
        \begin{tikzpicture}
            \draw[fill = lightgray,lightgray] (-1,-1)--(-1,1)--(1,1)--(1,-1)--(-1,-1);
            \draw[ultra thick, gray] (0,0) -- (0,1);
            \draw[ultra thick, gray] (0,0) -- (0,-1);
            \draw[ultra thick, gray] (0,0) -- (1,0);
            \draw[ultra thick, gray] (0,0) -- (-1,0);
            \draw[black,fill=black] (0,0) circle (.05cm);
            \draw (.6,.7) node {$NE$};
            \draw (-.6,.7) node {$NW$};
            \draw (-.6,-.7) node {$SE$};
            \draw (.6,-.7) node {$SW$};
            \draw[gray] (0,1.2) node {$N$};
            \draw[white,fill = white] (0,-1.38) circle (.1);
            \draw[gray] (0,-1.2) node {$S$};
            \draw[gray] (1.2,0) node {$E$};
            \draw[gray] (-1.2,0) node {$W$};
            \draw (.2,.2) node {$O$};
        \end{tikzpicture} \qquad
        \begin{tikzpicture}
            \draw[fill = lightgray,lightgray] (-1,-1)--(-1,1)--(1,1)--(1,-1)--(-1,-1);
            \draw[ultra thick, gray] (0,0) -- (0,1);
            \draw[ultra thick, gray] (0,0) -- (0,-1);
            \draw[ultra thick, gray] (0,0) -- (1,0);
            \draw[ultra thick, gray] (0,0) -- (-1,0);
            \draw[ultra thick, gray,fill=white] (0,0) circle (.4cm);
            \draw[black,fill=black] (0,.4) circle (.05cm);
            \draw[black,fill=black] (0,-.4) circle (.05cm);
            \draw[black,fill=black] (.4,0) circle (.05cm);
            \draw[black,fill=black] (-.4,0) circle (.05cm);
            \draw[white,fill = white] (0,-1.38) circle (.1);
        \end{tikzpicture} \qquad
        \begin{tikzpicture}
            \draw[white,fill=white] (0,-.95) circle (0.05cm);
            \draw[ultra thick, gray,fill=white] (0,0) circle (.4cm);
            \draw[black,fill=black] (0,.4) circle (.05cm);
            \draw[black,fill=black] (0,-.4) circle (.05cm);
            \draw[black,fill=black] (.4,0) circle (.05cm);
            \draw[black,fill=black] (-.4,0) circle (.05cm);
            \draw[white,fill = white] (0,-1.38) circle (.1);
        \end{tikzpicture}
    \end{center}
    \caption{Left: $(\RR^2 \to \cP)$ with indicated strata. Middle: the blowup $\bl_{\RR^2_O}(\RR^2 \to \cP)$. Right: the link $\link_{\RR^2_O}(\RR^2 \to \cP).$}
    \label{fig:plane}
\end{figure}

\begin{example}\label{cubeStrat}
    Let $\mc{P}$ be a poset. We define the set of \emph{subdivisions} of $\mc{P}$ to consist of all 
    subsets $S \subset \mc{P}$ such that $S$ is nonempty, finite, and the partial order on $S$ induced by 
    the partial order of $\mc{P}$ is a linear order. There is a natural partial order of the subdivisions of a poset
    given by inclusion of sets. We denote this poset by $\msf{sd}(\mc{P})$.
    The closed interval $[0,1]$ admits a natural stratification by the poset, $\msf{sd}(\{0<1\})$.
    The map $[0,1] \to \msf{sd}(\{0<1\})$ is given by
    $$x \mapsto \begin{cases} \{0\}, &\text{ if } x=0 \\ \{1\}, &\text{ if } x=1 \\ \{0,1\}, &\text{ else}
        \end{cases}~.$$
    More generally, let $k \in \ZZ_{>0}$. The following topological space
    $$\msf{Cube}^k := [0,1]^{\times k}~,$$
    admits a natural stratification by $\msf{sd}\left(\{0<1\}^{\times k}\right)$. This can be seen by taking the
    $k$-fold product of the above stratified space, $[0,1] \to \msf{sd}(\{0<1\})$, and then invoking the fact
    that $\left(\msf{sd}(\{0<1\})\right)^{\times k} \cong 
    \msf{sd}\left(\{0<1\}^{\times k}\right)$.
\end{example}

\subsection{Conical smoothness}
By imposing some extra regularity on a $\sC^0$-stratified manifold, we can use the stratification to compute the
homology of the manifold. This is the notion of a \emph{conically smooth structure}, which is analogous to a 
smooth structure on a topological space. 
We begin by briefly narrating key features of smooth structures on topological spaces (see~\cite{manifolds} for an introduction). \\

Among paracompact Hausdorff topological spaces, $\sC^0$-manifolds are characterized by being locally Euclidean, that is, each point has a neighborhood that is homeomorphic to $\RR^i$ for some $i\geq 0$.
A \emph{smooth structure} is a type of \emph{regularity} on a $\sC^0$-manifold.  
A $\sC^\infty$-manifold, or \emph{smooth} manifold, is a $\sC^0$-manifold equipped with a smooth structure.  
There is a distinguished class of continuous maps between two smooth manifolds, called the \emph{smooth} maps. This class of smooth maps consists precisely of those continuous maps that respect this smooth structure.
Here, ``respect'' is just so that points~(2) and~(3) are true.  
This regularity of a smooth manifold $M$ is tailored precisely so that it has the following features:
\begin{enumerate}
\item
For each point $x\in M$, there is a \emph{tangent space}, $T_x M$, which is a vector space.
This tangent space is a canonical local model of $M$ about $x$, which is to say there is a basis for the topology about $x\in M$ comprised of images of smooth open embeddings $\varphi_x: T_x M \hookrightarrow M$ each that carries $0$ to $x$  

\item
For $f: M \to N$ a smooth map, and for each $x\in M$ there is a linear map,
\[
D_x f: T_x M 
\longrightarrow 
T_{f(x)}N
~,
\]
called the \emph{derivative} of $f$ at $x$.  
Through a choice of smooth open embeddings $T_x M \xhookrightarrow{\varphi_x} M$ and $T_{f(x)}N \xhookrightarrow{\psi_{f(x)}} N$ as in the above point for which $f({\sf Image}(\varphi_x)) \subset {\sf Image}(\psi_{f(x)})$, this derivative can be identified as the limit
\begin{equation}\label{e1}
T_x M 
\xra{\varphi_x} 
{\sf Image}(\varphi_x) 
\xra{~f~}  
{\sf Image}( \psi_{f(x)})
T_{f(x)}N 
~,\qquad
v\mapsto 
\lim_{t\to 0} \frac{\psi_{f(x)}^{-1} f \varphi_x(tv)}{t}
~.
\end{equation}

\item
The map~(\ref{e1}) depends smoothly on $x\in M$, appropriately interpreted.  
\end{enumerate}

\begin{example}
Let $f=(f^1,\dots,f^n)\colon \RR^m \to \RR^n$ be a sequence of $n$ polynomial maps, each in $m$ variables.  
Suppose, for each $x\in f^{-1}(0)$, that the total Jacobian matrix of $f$ at $x$,
\[
\Bigl[ 
{\partial_i f^j}_{|x}
\Bigr]
~,
\]
has rank $n$.  
The \emph{Regular Value Theorem} states that the subspace $f^{-1}(0) \subset \RR^m$ has the natural structure of a smooth manifold.  
\end{example}

We now give a similar brief narrative of features of conically smooth structures on stratified spaces (see~\cite{aft1} for an account).
Among paracompact Hausdorff stratified topological spaces, $\sC^0$-stratified spaces are characterized by being locally a product of a Euclidean space and an open cone: $\RR^i\times \sC(L)$ for some $i\geq 0$ and compact $\sC^0$-stratified space.  
While this seems circular, because the topological dimension of $\sC(L)$ is strictly greater than that of $L$, this notion can be grounded through induction (on dimension).
A \emph{conically smooth} structure is a type of \emph{regularity} on a $\sC^0$-stratified space.  
A $\sC^\infty$-stratified space, or \emph{conically smooth} stratified space, is a $\sC^0$-stratified space equipped with a conically smooth structure.  
There is a special class of continuous maps between two conically smooth stratified spaces, called the \emph{conically smooth} maps. This class of maps consists precisely of those continuous maps which respect this conically smooth structure.
Here,``respect'' is just so that points~(2) and~(3) are true.  
This regularity of a conically smooth stratified space $X$ is tailored precisely so that it has the following features:
\begin{enumerate}
\item
For each point $x\in X$, there is a \emph{tangent cone}, $T_x X \times \sC(L_x)$, which is a product of a vector space and an open cone -- such a product has a scaling action by the group $\RR^\times$.  
This tangent cone is a canonical local model of $X$ about $x$, which is to say there is a basis for the topology about $x\in X$ comprised of images of smooth open embeddings $\varphi_x\colon T_x X \times \sC(L_x) \hookrightarrow X$ each that carries $0$ to $x$.

\item
For $f\colon X \to Y$ a conically smooth map, and for each $x\in X$ there is a $\RR^\times$-equivariant map,
\[
D_x f\colon T_x X \times \sC(L_x) 
\longrightarrow 
T_{f(x)} Y \times \sC(L_{f(x)})
~,
\]
called the \emph{derivative} of $f$ at $x$.  
Through a choice of smooth open embeddings $T_x X\times \sC(L_x) \xra{\varphi_x} X$ and $T_{f(x)}Y\times \sC(L_{f(x)}) \xra{\psi_{f(x)}} Y$ as in the above point for which $f({\sf Image}(\varphi_x)) \subset {\sf Image}(\psi_{f(x)})$, this derivative can be identified as the limit
\begin{equation}\label{e2}
T_x X \times \sC(L_x) 
\xra{\varphi_x} 
{\sf Image}(\varphi_x) 
\xra{~f~}  
{\sf Image}( \psi_{f(x)})
\xra{\psi_{f(x)}^{-1}} 
T_{f(x)} Y \times \sC(L_{f(x)})
~,
\end{equation}
\[
(v,[s,l])
~
\mapsto 
~
\lim_{t\to 0} \frac{\psi_{f(x)}^{-1} f \varphi_x(tv,[ts,l])}{t}
~.
\]

\item
The map~(\ref{e2}) depends conically smoothly on $x\in X$, appropriately interpreted.  
\end{enumerate}

\begin{example}
Let $f=(f^1,\dots,f^n)\colon \RR^m \to \RR^n$ be a sequence of $n$ polynomial maps each in $m$ variables.  
The Thom-Mather Theorem \cite{mather} states that, with no assumptions on the total Jacobian matrix, the subspace $f^{-1}(0) \subset \RR^m$ has the natural structure of a conically smooth stratified space.  
\end{example}

\begin{example}\label{r1}
Let $W\subset M$ be a properly embedded smooth $d$-submanifold of a smooth $n$-manifold.
Assume $d<n$.  
Consider the stratified topological space 
\[
(W\subset M) := \Bigl( M \xra{x\mapsto d ~\text{iff}~x\in W} \{d<n\} \Bigr)
~.
\]
So the $d$-stratum is $W = (W\subset M)_d$, and the $n$-stratum is $M\setminus W = (W\subset M)_n$.  
For each $x\in W$, consider the vector space $N_x:=\frac{T_x M}{T_x W}$.
Consider the unit sphere of this vector space $L_x := S(N_x) := ( N_x \setminus 0 )_{/\RR_{>0}}$ -- it is diffeomorphic with a $(n-d-1)$-sphere, and its open cone $\sC(L_x) \cong \RR^{n-d}$ is homeomorphic with $(n-d)$-Euclidean space.   
A choice of tubular neighborhood of $W\subset M$ determines, for each $x\in W$, an open embedding $\varphi_x \colon T_x W \times \sC(L_x) \hookrightarrow M$ that carries $0$ to $x$.  
These open embeddings determine a conically smooth structure on the stratified topological space $(W\subset M)$.  
\end{example}

The following is a catalogue of the key features of conically smooth stratified space $X = (X\to P)$:
\begin{enumerate}
\item
For each $p\in P$, the stratum $X_p$ is equipped with the structure of a connected smooth manifold.

\item
For each strictly related pair $p<q$ in $P$, there is a smooth manifold ${\sf Bl}_{X_p}(X)_q$ with boundary $\link_{X_p}(X)_q$ as well as a proper quotient map
\[
\ov{\pi}_{p<q} \colon {\sf Bl}_{X_p}(X)_q \longrightarrow X_p \cup X_q
\]
to the union in $X$ of the $p$- and the $q$-strata.  
This continuous map $\ov{\pi}_{p<q}$ has the following features.
\begin{enumerate}

\item
The preimage of the $p$-stratum is precisely the boundary of ${\sf Bl}_{X_p}(X)_q$; equivalently, the preimage of the $q$-stratum is precisely the interior of ${\sf Bl}_{X_p}(X)_q$:
\[
\ov{\pi}_{p<q}^{-1}(X_p) ~=~ \link_{X_p}(X)_q
\qquad\text{ and }\qquad
\ov{\pi}_{p<q}^{-1}(X_q) ~=~ {\sf Interior}\bigl( {\sf Bl}_{X_p}(X)_q\bigr)
~.
\]

\item
The restriction of $\ov{\pi}_{p<q}$ to the interior of ${\sf Bl}_{X_p}(X)_q$ is a diffeomorphism onto the $q$-stratum:
\[
(\ov{\pi}_{p<q})_{|} \colon {\sf Interior}\bigl({\sf Bl}_{X_p}(X)_q\bigr)
\xra{~\cong~}
X_q
~.
\]

\item
The restriction of $\ov{\pi}_{p<q}$ to the boundary of ${\sf Bl}_{X_p}(X)_q$ is a smooth fiber bundle:
\[
\pi_{p<q} \colon 
\link_{X_p}(X)_q
\longrightarrow
X_q
~.
\]

\end{enumerate}
Points (a)-(c) can be summarized as a commutative diagram
$$\begin{tikzcd}
    \link_{X_p}(X)_q \arrow[d,"\pi_{p<q}"] \arrow[rr,"\rm inclusion"]
    &&
    {\sf Bl}_{X_p}(X)_q \arrow[d,"{\ov{\pi}}_{p<q}"]
    &&
    {\sf Interior} \bigl( {\sf Bl}_{X_p}(X)_q \bigr) \arrow[ll,swap,"\rm inclusion"] \arrow[d,"\cong"]
    \\
    X_p \arrow[rr,"\rm inclusion"]
    &&
    X_p \cup X_q
    &&
    X_q \arrow[ll,swap,"\rm inclusion"]
\end{tikzcd}$$
in which each square is both a pullback and the left square is a pushout.  
In particular, a choice of collaring of the boundary,
\[
\ov{\gamma}_{p<q}\colon \link_{X_p}(X)_q \times [0,1)
\hookrightarrow
{\sf Bl}_{X_p}(X)_q
~,
\]
restricts as a smooth open embedding
\begin{equation}\label{e3}
\gamma_{p<q}\colon 
\link_{X_p}(X)_q \times (0,1)
\hookrightarrow
X_q
~.
\end{equation}

\item
For strictly related triples $p<q<r$ in $P$, there is a smooth manifold with $\lag 2 \rag$-corners ${\sf Bl}_{X_p\cup X_q}(X)_r$, together with a proper quotient map $\ov{\pi}_{p<q<r}\colon {\sf Bl}_{X_p\cup X_q}(X)_r\to X_p\cup X_q \cup X_r$ with similar features to (a)-(c) above.  

\item
Etcetera, for finite strictly monotonic sequences $p_1<\dots<p_\ell$ in $P$. 
\end{enumerate}

\begin{example}\label{r2}
We follow up on Example~\ref{r1}.
Namely, the smooth structures on the strata $W = (W\subset M)_d$ and on $M\setminus W = (W\subset M)_n$ are the given ones inherited from the given smooth structure on $M$.  
The link 
\[
\Bigl(
\link_{W}\bigl( (W\subset M) \bigr)  \xra{~\pi_{d<n}~} W
\Bigr)
~=~
\Bigl(
S^{\sf fib}(N_{W\subset M}) \xra{~\pr~} W
\Bigr)
\]
is identical with the unit sphere bundle of the normal bundle of $W\subset M$, which is a smooth manifold.
The smooth $n$-manifold with boundary ${\sf Bl}_W\bigl( (W\subset M) \bigr)$ is the \emph{real} blow-up of $M$ along $W$.  
The interior of this real blow-up is the complement $M\setminus W$.  
See Figure~\ref{fig:R2} in the case that $(W\subset M) = (\{0\}\subset \RR^2)$. 
\end{example}

\begin{figure}
    \begin{center}
    \begin{tikzpicture}
    \draw[lightgray, fill=lightgray] (-1,1)--(-1,-1)--(1,-1)--(1,1)--(-1,1);
    \draw[black,fill=black] (0,0) circle (.075);
    \end{tikzpicture}\qquad
     \begin{tikzpicture}
    \draw[lightgray, fill=lightgray] (-1,1)--(-1,-1)--(1,-1)--(1,1)--(-1,1);
    \draw[ultra thick,black,fill=white] (0,0) circle (.5);
    \end{tikzpicture}\qquad
    \begin{tikzpicture}
    \draw[white, fill=white] (-1,1)--(-1,-1)--(1,-1)--(1,1)--(-1,1);
    \draw[ultra thick,black,fill=white] (0,0) circle (.5);
    \end{tikzpicture}
    \end{center}
    \caption{Left: the stratified space $(\{0\} \subset \RR^2)$. Middle: the blowup $\bl_{\{0\}}\bigl( (\{0\}\subset \RR^2) \bigr)$. Right: the link $\link_{\{0\}}\bigl( (\{0\} \subset \RR^2) \bigr)$.}
    \label{fig:R2}
\end{figure}

\begin{notation}
    Let $X \to \mc{P}$ be a conically smooth stratified space, and let $p<q$ be a pair of strictly related elements of
    $\mc{P}$. We will often denote the link by
    $$\msf{L}_{p<q}(X) := \link_{X_p}(X)_q~.$$
\end{notation}

\subsection{Links and homology}
We now explain how the previously defined links allow one to compute the homology of a conically smooth 
stratified space.
In this section, we restrict attention to topological spaces, $X$, that possess the following structure:
\begin{itemize}
    \item $X$ is a smooth manifold;
    \item $X$ is equipped with a conically smooth structure $X \to \mc{P}$ that is compatible with the smooth
        structure;
    \item Each stratum, $X_p$, is diffeomorphic to Euclidean space, and we further fix such a diffeomorphism 
        $$\alpha_p:X_p \xrightarrow{\cong} \bb{R}^{\msf{dim}(X_p)}~,$$
        for each $p \in \mc{P}$.
\end{itemize}
This structure gives rise to a well-defined map of posets
$$\msf{d}:\mc{P} \to \bb{Z}_{\ge 0}~, \qquad p \mapsto \dim(X_p)~.$$
We let $X_{(i)}$ denote the fiber of the composite $X \to \mc{P} \to \bb{Z}_{\ge 0}$:
$$\begin{tikzcd}
        X_{(i)} \arrow[r,hook] \arrow[dd] \arrow[dr, phantom, "\lrcorner", very near start] & X \arrow[d] \\
        & P \arrow[d] \\
        \{i\} \arrow[r,hook] & \bb{Z}_{\ge 0}
    \end{tikzcd}~,$$
and let $X_{(\le i)}$ denote the pullback
$$\begin{tikzcd}
        X_{(\le i)} \arrow[r,hook] \arrow[dd] \arrow[dr, phantom, "\lrcorner", very near start] & X \arrow[d] \\
        & P \arrow[d] \\
        \bb{Z}_{\le i} \arrow[r,hook] & \bb{Z}_{\ge 0}
    \end{tikzcd}~.$$

Fixing diffeomorphisms, $\alpha_p$ for each $p \in \mc{P}$, endows each
stratum with an orientation, and thus an isomorphism $\redc_*(X_p^+) \cong_{\alpha_p} \bb{Z}[\dim(X_p)]$ between the reduced 
chain complex of $X_p^+$ and $\bb{Z}$ in dimension $\dim(X_p)$. So to each such 
conically smooth stratified space, we obtain a 
sequence of abelian groups
$$\begin{tikzcd}
    \bb{Z}\langle \pi_0(X_{(0)}) \rangle [0] & \bb{Z}\langle \pi_0(X_{(1)}) \rangle [1] & \cdots & \bb{Z}\langle \pi_0(X_{(d)})
        \rangle [d] & 0 & \cdots
    \end{tikzcd}~,$$
where $\bb{Z}\langle \pi_0(X_{(i)})\rangle$ sits in degree $i$. \\

To name
a homomorphism
$$\bb{Z}\langle \pi_0(X_{(i+1)})\rangle \xrightarrow{\partial_{i+1}} \bb{Z}\langle \pi_0(X_{(i)})\rangle$$
is to specify a $\pi_0(X_{(i)}) \times \pi_0(X_{(i+1)})$-matrix. Each element of 
$\pi_0(X_{(i)})$ is a stratum, $X_p$, of $X$ with $\dim(X_p) = i$.
First, for each $p<q$ in $\mc{P}$, such that $\dim(X_p) +1=\dim(X_q)$, we define a map
$$\sigma_p^q: \pi_0(\msf{L}_{p<q}(X)) \to \{\pm1\}$$
by
$$[\ell] \mapsto \msf{det}\left(\bb{R}^i \oplus \bb{R} \cong T_{\pi(\ell)} X_p \oplus \bb{R} 
    \xrightarrow{(D_\ell \pi)^{-1} \oplus \text{id}} T_\ell \msf{L}_{p<q}(X) \oplus \bb{R} \cong
    T_{(\ell,\varepsilon/2)} (\msf{L}_{p<q}(X) \times (0,\varepsilon)) \xrightarrow{D_\ell \gamma} T_{\gamma(\ell)}X_q 
    \cong \bb{R}^{i+1}\right)~.$$
Now, we define the $(p,q)$ entry of 
$\partial_{i+1}:\bb{Z}\langle \pi_0(X_{(i+1)})\rangle \to \bb{Z}\langle \pi_0(X_{(i)})\rangle$ to be 
\begin{equation}\label{differential}
    \sum_{[\ell] \in \pi_0(\msf{L}_{p<q}(X))} \sigma_p^q([\ell])~.
\end{equation}

\begin{proposition}\label{homology}
    Let $X \to \mc{P}$ be a compact, connected, conically smooth stratified manifold
    such that each stratum is diffeomorphic to a Euclidean space. The previously defined sequence of abelian
    groups, and homomorphisms between them
    $$\partial_{i+1}: \bb{Z}\langle \pi_0(X_{(i+1)})\rangle \to \bb{Z}\langle \pi_0(X_{(i)})\rangle~,$$
    is a chain complex. Furthermore, the homology of this chain complex is isomorphic with 
    $H_*(X;\bb{Z})$, the singular homology of $X$ with $\bb{Z}$ coefficients.
\end{proposition}
{\it Proof.} The filtration
$$X_{(\le0)} \hookrightarrow X_{(\le 1)} \hookrightarrow \cdots \hookrightarrow X_{(\le \ell)} \hookrightarrow
    \cdots \hookrightarrow X~,$$
of $X$ induces a filtration of $C_*(X)$,
$$C_*(X_{(\le0)}) \to C_*(X_{(\le 1)}) \to \cdots \to C_*(X_{(\le \ell)}) \to \cdots \to C_*(X)~.$$
The $E^0$ page of the spectral sequence associated to this filtration is
$$E^0_{i,j} = C_{i+j}(X_{(\le i)}) / C_{i+j}(X_{(\le i-1)})~.$$
Hence the $E^1$ page is given by
$$E^1_{i,j} = H^\text{red}_{i+j}\left(C_{i+j}(X_{(\le i)}) / C_{i+j}(X_{(\le i-1)}\right) \cong
    H^\text{red}_{i+j}\left(X_{(\le i)} / X_{(\le i-1)}\right)~,$$
where the homology is taken with respect to the $d^1$-differential.
Since $X$ is compact and connected, $X_{(\le i)} / X_{(\le i-1)} \cong (X_i)^+$, where
$(X_{(i)})^+$ denotes the one-point compactification of $X_{(i)}$. Further, since each stratum is equipped with a diffeomorphism $\alpha_p:X_p \xrightarrow{\cong} \bb{R}^{\msf{dim}(X_p)}$, 
we have the based homeomorphism,
\begin{equation}\label{wedge}
    X_{(i)}^+ \cong_{\alpha_p^+} \left(\coprod_{\pi_0(X_{(i)})} \RR^i \right)^+
    \cong \bigvee_{\pi_0(X_{(i)})} S^i~.
\end{equation}
Therefore, 
$$E^1_{i,j} \cong H^\text{red}_{i+j}\left( \bigvee_{\pi_0(X_{(i)})} S^i\right)
    \cong \begin{cases} \bb{Z}\langle \pi_0(X_{(i)})\rangle , &\quad \text{if } j=0 \\
                                      0, &\quad \text{otherwise}
              \end{cases}~.$$
The $d^1$ differential is then a homomorphism
$$d^1_{i,0}: \bb{Z}\langle \pi_0(X_{(i)})\rangle \to \bb{Z}\langle \pi_0(X_{(i-1)})\rangle~.$$
We will now identify $d^1$ with the homomorphism (\ref{differential}). This will complete the proof of the first 
statement, since the spectral sequence collapses at the $E^2$ page for dimension reasons. 
Further, convergence of the spectral sequence
to $H_*(X; \bb{Z})$ is seen from the fact that the filtration is finite. \\

The $d^1$ differential is given by the map induced by applying $H_*^\text{red}$ to the following composite morphism of spaces
$$\frac{X_{(\le i)}}{X_{(\le i-1)}} \simeq \cone\left(X_{(\le i-1)} \hookrightarrow X_{(\le i)}\right) \to 
    \Sigma X_{(\le i-1)}^+ \to \Sigma \frac{X_{(\le i-1)}}{X_{(\le i-2)}}~,$$
where the first morphism is given by collapsing $X_{(\le i)}$ in the cone, and the second morphism is the quotient.
Consider the continuous map
$$\gamma_{(\le i-1)}^!: X_{(i)}^+ \to \Sigma \link_{X_{(\le i-1)}}(X_{(\le i)})^+~,$$
which is the following composite
$$\begin{tikzcd}
    X^+_{(i)} \arrow[d,"\cong",swap] \arrow[rrr,"\gamma_{(\le i-1)}^!",dashed] & & & 
        \Sigma \link_{X_{(\le i-1)}}(X_{(\le i)})^+ \\
    \frac{X_{(\le i)}}{X_{(\le i-1)}} \arrow[r,"\cong"] 
        & \frac{\bl_{X_{(\le i-1)}}(X_{(\le i)})}{\link_{X_{(\le i-1)}}(X_{(\le i)})} \arrow[r,"\text{collapse}"]
        & \frac{(\link_{X_{(\le i-1)}}(X_{(\le i)}) \times [0,1))^+}{\link_{X_{(\le i-1)}}(X_{(\le i)}) \times \{0\}}
            \arrow[r,"\cong"]
        & (\link_{X_{(\le i-1)}}(X_{(\le i)}) \times (0,1))^+ \arrow[u,"\cong"]
    \end{tikzcd}~.$$
Consider the following solid commutative diagram
$$\begin{tikzcd}
    X_{(i)}^+ \arrow[r,"\gamma^!_{(\le i-1)}"] \arrow[dr] \arrow[rr,bend left=40,"\gamma^!_{(i-1)}",dashed]
        & \Sigma \link_{X_{(\le i-1)}}(X_{(\le i)})^+ \arrow[d,"\Sigma\pi_{(\le i-1)}"] \arrow[r,"\text{collapse}"]  
        & \Sigma\link_{X_{(i-1)}}(X_{(i-1)}\cup X_{(i)})^+ \arrow[d,"\Sigma\pi_{(i-1)}"] \\
    & \Sigma X_{(\le i-1)}^+ \arrow[r] & \Sigma X_{(i-1)}^+
    \end{tikzcd}~,$$
Denote the indicated horizontal composite as $\gamma^!_{(i-1)}$. Since the bottom
composite induces the $d^1$ differential upon applying reduced homology, we have identified the $d^1$
differential with the map induced on homology by the composite $\Sigma\pi_{(i-1)} \circ \gamma^!_{(i-1)}.$
Using the fact that for $Y$ a locally compact and Hausdorff space, there is a based homeomorphism 
$\Sigma Y^+ \cong (Y \times (0,1))^+$, we have the based homeomorphism
\begin{align*}
    \Sigma\link_{X_{(i-1)}}(X_{(i-1)} \cup X_{(i)})^+ 
    &\cong \left(\link_{X_{(i-1)}}(X_{(i-1)} \cup X_{(i)}) \times (0,1)\right)^+ \\
    &\cong \left( \coprod_{\pi_0(\link_{X_{(i-1)}}(X_{(i-1)} \cup X_{(i)}))} \RR^{i-1} \times (0,1)\right)^+ \\
    &\cong \bigvee_{\pi_0(\link_{X_{(i-1)}}(X_{(i-1)} \cup X_{(i)}))} S^i~.
\end{align*}
Similar to the above identification we also have the following based homeomorphism
$$\Sigma X_{(i-1)}^+
    \cong \bigvee_{\pi_0(X_{(i-1)})} S^i~.$$
Using these identifications together with Equation (\ref{wedge}), we identify 
$\gamma^!_{(i-1)}$  with a map between wedges of spheres
$$\bigvee_{\pi_0(X_{(i)})} S^i \cong_{\alpha_i^+} X_{(i)}^+ \xrightarrow{\gamma^!_{(i-1)}}
    \Sigma\link_{X_{(i-1)}}(X_{(i-1)} \cup X_{(i)})^+\xrightarrow{\cong}
    \bigvee_{\pi_0(\link_{X_{(i-1)}}(X_{(i-1)} \cup X_{(i)}))} S^i~,$$
and we identify $\Sigma\pi_{(i-1)}$ as a map between wedges of spheres, 
$$\bigvee_{\pi_0(\link_{X_{(i-1)}}(X_{(i-1)} \cup X_{(i)}))} S^i \cong 
    \Sigma\link_{X_{(i-1)}}(X_{(i-1)} \cup X_{(i)})^+ \xrightarrow{\Sigma\pi_{(i-1)}} \Sigma X_{(i-1)}^+
    \cong \bigvee_{\pi_0(X_{(i-1)})} S^i~.$$
Thus, upon applying $H^\text{red}_*$, the $d^1$ differential is the composite map
$$\bb{Z}\langle \pi_0(X_{(i)})\rangle \to 
    \bb{Z}\langle \pi_0(\link_{X_{(i-1)}}(X_{(i-1)}\cup X_{(i)})^+)\rangle \to 
    \bb{Z}\langle \pi_0(X_{(i-1)})\rangle~.$$
We can compute these induced maps via the degrees of the maps $\gamma^!_{(i-1)}$ and $\Sigma\pi_{(i-1)}$
(see \cite{hatcher} for more details). To identify the $d^1$-differential with Equation (\ref{differential}),
and thus complete the proof, we now identify the degree of 
$$\Sigma\pi_{(i-1)} \circ \gamma^!_{(i-1)} : \bigvee_{\pi_0(X_{(i)})} S^i \to
    \bigvee_{\pi_0(X_{(i-1)})} S^i~,$$
with the formula in Equation (\ref{differential}). 
The degree of a smooth map $f:X \to Y$ between compact, oriented smooth manifolds can be computed by 
choosing a regular value $y \in Y$, and then taking a count of preimages of $y$, signed according to whether 
$f$ is orientation preserving or reversing at that point:
$$\msf{deg}(f) = \sum_{x \in f^{-1}(y)} \msf{sgn}\left(\msf{det}(D_x f)\right)~,$$
(see \cite{guillemin}, for instance). 
For ease of notation, 
let us denote the composite $\Sigma\pi_{(i-1)} \circ \gamma^!_{(i-1)}$  by $\Phi$. Note that $\Phi$ is smooth
away from the basepoint. Further, each element of $\pi_0(X_{(i)})$ is a stratum $X_q$, for some $q \in \mc{P}$,
such that $d(X_q) = i$. Similarly, each element of $\pi_0(X_{(i-1)})$ is a stratum $X_p$, for some $p \in \mc{P}$,
such that $d(X_p) = i-1$. Thus, we can compute the degree by computing the degree of the induced map
$$\Phi_{p<q}:S^i_q := X_q^+ \to (X_p \times (0,1))^+=:S^i_p~,$$
for each $p<q \in \cP$ with $d(p)+1=d(q)$. Choosing a regular value, $y \in S^i_p$, we have
$$\msf{deg}(\Phi_{p<q}) = \sum_{x \in \Phi_{p<q}^{-1}(y)} 
    \msf{sgn}\left(\msf{det}(D_x \Phi_{p<q})\right)~.$$
Observing that $\card(\Phi_{p<q}^{-1}(y)) = \card(\pi_0(\msf{L}_{p<q}(X)))$, it just remains to establish that 
$$\msf{sgn}(\msf{det}(D_x \Phi_{p<q})) = \msf{sgn}(\msf{det}(D_\ell \gamma \circ ((D_\ell \pi)^{-1} \oplus 
    \msf{id}_{(0,\varepsilon)})))~,$$
as in Equation (\ref{differential}). Since $\Phi_{p<q}$ is given as the composite, we have
$$D_x \Phi_{p<q} = (D_\ell \pi \oplus \msf{id}_{(0,\varepsilon)}) \circ (D_x \gamma)^{-1}~,$$
where $\ell = \gamma^{-1}_{p<q}(x)$. Hence the sign of this determinant matches that of Equation 
(\ref{differential}), which completes the proof.
\qed \\

\section{The Schubert stratification}\label{sec:schubert}
In Section \ref{cell}, we showed that the cells of the Grassmannian $\gr_k(n)$ are parameterized by the set
$$\nchoosek{n}{k}:=\{S \subset \{1,\dots,n\} ~|~ \card(S) = k\}~.$$
Since each element in $\nchoosek{n}{k}$ corresponds to a cell in the CW structure on $\gr_k(n)$, there is a 
canonical map 
$$\gr_k(n) \to \nchoosek{n}{k}$$
given by sending a $k$-plane to the set indexing its
cell. This realizes $\gr_k(n)$ as a stratfied topological space with $S$-stratum the $S$-cell:
$$\gr_k(n)_S := \{ V \in \gr_k(n) ~|~ 
    \text{$S$ is the maximal element in } {\nchoosek{n}{k}} \text{ for which }V \in U_S \} ~.$$
In particular, recall that each stratum of $\gr_k(n)$ is diffeomorphic to Euclidean space. 
In fact, as proven in a forthcoming paper, we prove:
\begin{theorem}\label{conicallySmooth}
    For $0 \le k \le n$, the Schubert stratification of the Grassmannian $\gr_k(n)$ can be naturally upgraded to
    the structure of a conically smooth stratified space, with respect to which, for each pair $S<T$ in 
    $\nchoosek{n}{k}$ with $d(S)+1=d(T)$, there are identifications
    $$\begin{tikzcd}
        \gr_k(n)_S & \msf{L}_{S<T}(\gr_k(n)) \arrow[l,"\pi_{S<T}"'] \\
        (0,\pi)^{Z_S} \arrow[u,"\cong","\rot_S"'] & \msf{L}_{S<T} \arrow[l,"\msf{proj}"'] \arrow[u,"\cong"]
    \end{tikzcd} \text{ and }
    \begin{tikzcd}
        \msf{L}_{S<T}(\gr_k(n)) \times (0,\pi) \arrow[r,"\gamma_{S<T}"] & \gr_k(n)_T \\
        \msf{L}_{S<T} \times (0,\pi) \arrow[u,"\cong"] \arrow[r,"\msf{swap} \amalg \msf{Rswap}"] 
            & (0,\pi)^{Z_T} \arrow[u,"\cong","\rot_T"']
    \end{tikzcd} ~,$$
    where 
    \begin{align*}
        \msf{L}_{S<T} &:= \left\{ \left(\theta_{(i,j)}\right) \in [0,\pi]^{Z_T} \big| \text{ for all } (i,j) \in Z_T, 
        \text{ if } (i,j) \in Z_S, \theta_{(i,j)} \in (0,\pi), \text{ else } \theta_{(i,j)} \in \{0,\pi\} \right\} \\
        &\cong (0,\pi)^{Z_S} \times \{0,\pi\}~.
    \end{align*}
\end{theorem}


As stated by Proposition \ref{homology}, to compute the homology of $\gr_k(n)$, we need to understand
the differential given in Equation (\ref{differential}). In particular, for 
$S<T \in \nchoosek{n}{k}$ with $d(S) +1 = d(T)$, we must understand the map
$$\sigma_S^T: \pi_0(\msf{L}_{S<T}) \to \{\pm1\}$$
defined by sending $[\ell] \in \pi_0(\msf{L}_{S<T})$ to 
$$\msf{det}\left(\bb{R}^i \oplus \bb{R} \cong T_{\pi_S^T(\ell)} \gr_k(n)_S \oplus \bb{R} 
    \xrightarrow{(D_\ell \pi_S^T)^{-1} \oplus \text{id}} T_\ell \msf{L}_{S<T} \oplus \bb{R} \cong
    T_{(\ell,\varepsilon/2)} (\msf{L}_{S<T} \times (0,\varepsilon)) \xrightarrow{D_\ell \gamma_S^T} T_{\gamma_S^T(\ell)}\gr_k(n)_T \cong \bb{R}^{i+1}\right)$$
$$=\msf{det}\left((D_\ell \pi_S^T)^{-1}\right) \cdot \msf{det}\left(D_\ell \gamma_S^T\right)~.$$
So, let us fix $S<T \in \nchoosek{n}{k}$ with $d(S)+1=d(T)$.  \\

The following immediate consequence of Corollary \ref{conicallySmooth} tells us that we only need to compute 
$\sigma^T_S$ for two values.

\begin{corollary}\label{piNaughtLink}
    For $S<T \in \nchoosek{n}{k}$ with $d(S)+1=d(T)$, $\pi_0(\msf{L}_{S<T}) \cong \{0, \pi\}$.
\end{corollary}

Let us fix two such values, one in each connected component of $\msf{L}_{S<T}$:
$$\Theta_{S<T}^0 = (\theta^0_{(i,j)})_{(i,j)\in Z_T}~, \qquad \theta^0_{(i,j)} := 
    \begin{cases} \pi/2, &\text{ if } (i,j) \in Z_S \\ 0, &\text{ if } (i,j) \in Z_T\setminus Z_S\end{cases}~,$$
and 
$$\Theta_{S<T}^\pi = (\theta^\pi_{(i,j)})_{(i,j)\in Z_T}~, \qquad \theta^\pi_{(i,j)} := 
    \begin{cases} \pi/2, &\text{ if } (i,j) \in Z_S \\ \pi, &\text{ if } (i,j) \in Z_T\setminus Z_S \end{cases}~.$$

We will now unpack the maps $D_\ell \pi_S^T$ and $D_\ell\gamma_S^T$. Let us now fix $S<T \in \nchoosek{n}{k}$
with $d(S)+1=d(T)$.
Let $M \in \nchoosek{n}{k}$ be the
maximal element. \\

We have the following commutative diagrams of spaces:
\begin{center}
    \begin{tikzcd}
        \msf{L}_{S<T}\arrow[rr,"\pi_{S<T}"]\arrow[d,"\msf{pr}"]&&\gr_k(n)_S
            \arrow[d,hook,"\inc_S"]\\
        \left(0,\pi\right)^{Z_S}\arrow[rr,"\rot_S"]\arrow[dr,"\trot_S"',sloped]&&\gr_k(n)\\
        &\msf{V}_k^o(n)\arrow[ur,"\col"',sloped]&
    \end{tikzcd}~,
    \begin{tikzcd}
        \msf{L}_{S<T}\times(0,\pi)\arrow[rr,"\gamma_{S<T}"]\arrow[d,"\msf{swap}\amalg\msf{Rswap}"]
            &&\gr_k(n)_T \arrow[d,hook,"\inc_S"]\\
        \left(0,\pi\right)^{Z_T}\arrow[rr,"\rot_T"]\arrow[dr,"\trot_T"',sloped]&&\gr_k(n)\\
        &\msf{V}_k^o(n)\arrow[ur,"\col"',sloped]&
    \end{tikzcd}~,
\end{center}

from which, we compute $D_\ell \pi_{S<T}$ and $D_\ell \gamma_{S<T}$ as their respective composites.

\begin{lemma}\label{boundarySchubert}
    For $S<T \in \nchoosek{n}{k}$ with $d(S)+1=d(T)$, 
    $$\sigma_S^T([\ell]) = \begin{cases}
                                            (-1)^{1+k+i_{S<T}+\sum\limits_{i=i_{S<T}+1}^k s_i-i}, &\quad \text{if } \ell \in [0] \in 
                                                \pi_0(\msf{L}_{S<T}) \\
                                            (-1)^{\sum\limits_{i=i_{S<T}}^k s_i-i}, &\quad
                                                \text{if } \ell \in [\pi] \in \pi_0(\msf{L}_{S<T})
                                        \end{cases}~,$$
    where $1 \le i_{S<T} \le k$ is the index for which $S_{i_{S<T}}+1=T_{i_{S<T}}$.
\end{lemma}
{\it Proof.}
Let us first consider the derivative of $\col$. 
For $A \in V_k^o(n)$, $T_AV_k^o(n) = \{V \in \mat_{n\times k} ~|~ V^TA + A^TV = 0\}$,
and $T_{\col(A)}(\gr_k(n)) \cong \hom(\col(A), \col(A)^\perp) \cong \mat_{(n-k) \times k}$. Thus, a basis
for $T_{\col(A)}(\gr_k(n))$ is given by the set $\mc{M} := \{\msf{M}_{(i,j)} ~|~ 1 \le i \le n \setminus k, ~1 \le j
\le k\}$, where $\msf{M}_{(i,j)}$ has a 1 in the $(i,j)$ entry and all other entries are 0.
Thus, for $\col(A) \in \gr_k(n)_S$, $T_{\col(A)}(\gr_k(n)_S)$ has a basis given by
$$\mc{M}_S := \{M_{(i,j)} ~|~ 1 \le j \le k, ~ 1 \le i < s_j\} \subset \mc{M}.$$
Note that we will refer to elements of $\mc{M}$ or $\mc{M}_S$ as the indexing pairs, $(i,j)$.Recall that $\col:V_k^o(n \to \gr_k(n)$ takes an $n\times k$
matrix to its column space. Let $A \in \{\trot_S(\Theta_{S<T}^0),\trot_S(\Theta_{S<T}^\pi)\}$.
The surjective submersion $O(n) \to V_k^o(n)$ tells us that for each 
$V \in T_AV_k^o(n)$, there
exists $\widetilde{V} \in \msf{Skew}(n)$ for which $e^{t\widetilde{V}}A$ represents the equivalence class of the
tangent vector $V$. 
Let $\uno_{(\ul{n}\setminus S) \times n}$ denote the $n\times n$ identity matrix with the rows labeled by $S$ removed.
Now, in local coordinates, we can express $\col$ as the matrix 
$$\uno_{(\ul{n}\setminus S) \times n} e^{t\widetilde{V}}A(A^Te^{t\widetilde{V}}A)^{-1}~,$$
whose derivative at $t=0$ is
$$\frac{d}{dt}\bigg|_{t=0} \left(\uno_{(\ul{n}\setminus S) \times n} e^{t\widetilde{V}}A(A^Te^{t\widetilde{V}}A)^{-1}
    \right) = \uno_{(\ul{n}\setminus S) \times n} \widetilde{V} A~.$$
When $A = \trot_S(\Theta_{S<T}^0)$, $\widetilde{V}A = V$, so $D_{\Theta_{S<T}^0} \col : V_k^o(n) \to \gr_k(n)$
sends $V$ to the $\ul{n}\setminus S \times k$ matrix that consists of the $\ul{n}\setminus S$ rows of $V$. When
$A = \trot_S(\Theta_{S<T}^\pi)$, $\widetilde{V}A$ is $V$ except with the $i_{S<T}$ column negated.
Thus, $D_{\Theta_{S<T}^\pi} \col$ sends $V$ to the matrix that consists of the $\ul{n}\setminus S$ rows of $V$, but
with the $i_{S<T}$ column negated. Note that 
$T_AV_k^o(n) = \{V \in \msf{Mat}_{n\times k} ~|~ V^TA+A^TV=0\}$ has a natural basis consisting of two types
of matrices. The first type consists of matrices whose only nonzero entry is the $(i,j)$ entry, for 
$(i,j) \not\in S \times S$, one such matrix for each such $(i,j)$. 
The second type of matrices will have their nonzero entries concentrated in the rows and columns labeled by $S$. 
Note that our description of the derivative shows that for both $0$ and $\pi$, the derivative sends all matrices of
the second type to 0. Note that $T_{\col(A)}\gr_k(n) \cong \msf{Mat}_{(n-k) \times k}$ has a natural
basis consisting of the matrices of the first type described above. Thus, deleting the second type of basis elements
of $T_AV_k^o(n)$, $D_{\trot_S(\Theta_{S<T}^0)}\col$ is simply the identity matrix, and 
$D_{\trot_S(\Theta_{S<T}^\pi)}$ is the identity matrix except with all columns labeled by $(i,j)$ for which 
$j=i_{S<T}$. Finally, since we are only interested in those basis elements for which $(i,j) \in \mc{M}_{S}$, 
as discussed above, we conclude that
$$\det(D_{\trot_S(\Theta_{S<T}^0)}\col) = 1~,$$
and
$$\det(D_{\trot_S(\Theta_{S<T}^\pi)}\col) = (-1)^{\card\{(i,j) \in \mc{M}_{S^c}^* ~|~ j=i_{S<T}\}}
    = (-1)^{s_{i_{S<T}}-i_{S<T}}~.$$

Note that 
$\trot_S(\Theta^0_{S<T}) =\trot_S(\Theta^\pi_{S<T}) \in V^o_k(n)$ is the $n \times k$ matrix whose $i$th 
column is $\e_{s_i}$. 
We will first consider
$$D_{\Theta^0_{S<T}} \trot_S = \left[ \left(\frac{\partial{\trot_S}}{\partial \theta_{(i,j)}}\right)_{xy} \right]~.$$
That is, for each pair, $(i,j) \in Z_S$ and $(x,y) \in \mc{M}_S$, we will identify the $xy$-entry in the 
matrix $\frac{\partial{\trot_S}}{\partial \theta_{(i,j)}}$. Note that we ignore those entries indexed by elements
$(x,y) \in \mc{M} \setminus \mc{M}_S$, as $D\col$ will map those to 0, as discussed above.
Let us now fix $(i,j) \in Z_S$. Then
$$\frac{\partial{\trot_S}}{\partial \theta_{(i,j)}}(\Theta_{S<T}^0) = 
    \left( \prod_{(i',j') < (i,j)} \msf{R}_{j'}(\theta_{(i',j')})\right) D\msf{R}_j(\theta_{(i,j)}) 
    \left( \prod_{(i'',j'')>(i,j)} \msf{R}_{j''}(\theta_{(i'',j'')})\right)~.$$
For ease of notation, let us denote
$$\msf{Rot}^{<(i,j)}_S := \prod_{(i',j') < (i,j)} \msf{R}_{j'}(\theta_{(i',j')})~,$$
and 
$$\msf{Rot}^{>(i,j)}_S:=\prod_{(i'',j'')>(i,j)} \msf{R}_{j''}(\theta_{(i'',j'')})~.$$
Recall that $\theta_{(i',j')} = \pi/2$ for each $(i',j') \in Z_S$, so each entry in $D\msf{R}_j(\theta_{(i,j)})$ is zero,
except for the $(j,j)$ and $(j+1,j+1)$ entries which are both $-1$. Thus, we only need to know which $\e_r$ gets
sent to $\e_j$ and $\e_{j+1}$ under the map $\msf{Rot}^{>(i,j)}_S$. By definition, $\msf{Rot}^{>(i,j)}_S$
is a product of cyclic permutation matrices
$$[\e_i \mapsto \e_j] [\e_{i+1} \mapsto \e_{s_{i+1}}] \cdots [\e_k \mapsto \e_{s_k}]~.$$
The first block $[\e_i \mapsto \e_j]$ sends
$$\e_r \mapsto \begin{cases} \e_j, &\text{ if } r=i \\
                                                  \e_{j+1}, &\text{ if } r=j+1
                         \end{cases} ~.$$
For the next block, $[\e_{i+1} \mapsto \e_{s_{i+1}}]$, we care about what gets sent to $\e_i$ and $\e_{j+1}$.
We see that $[\e_{i+1} \mapsto \e_{s_{i+1}}]$ sends
$$\e_r \mapsto \begin{cases} \e_i, &\text{ if } r=i \\
                                                  -\e_{j+1}, &\text{ if } r=j+2 
                         \end{cases} ~.$$
Similarly, the next block will fix $\e_i$, and send $\e_{j+3}$ to $-\e_{j+2}$. Since there are $k-i$ blocks after
$[\e_i \mapsto \e_j]$, we see that
$$D\msf{R}_j \circ \msf{Rot}^{>(i,j)}_S : \e_r \mapsto \begin{cases} -\e_j, &\text{ if } r=i \\
                                                                                          (-1)^{k-i+1}\e_{j+1}, &\text{ if } r=j+k-i+1 \\
                                                                                          0, &\text{ else}
                                                                    \end{cases} ~.$$
Next, we will determine where $\msf{Rot}^{<(i,j)}_S$ sends $\e_j$ and $\e_{j+1}$. Note that
$\widehat{\rot}_S$ is the following product of cyclic permutation matrices
$$[\e_1 \mapsto \e_{s_1}] \cdots [\e_{i-1} \mapsto \e_{s_{i-1}}] [\e_{j+1} \mapsto \e_{s_i}].$$
The $[\e_{j+1} \mapsto \e_{s_i}]$ block fixes $\e_j$ and sends $\e_{j+1} \mapsto \e_{s_i}$. The next block from
the right, $[\e_{i-1} \mapsto \e_{s_{i-1}}]$ sends
$$\e_j \mapsto \begin{cases} \e_j, &\text{ if } j>s_{i-1} \\
                                                -\e_{j-1}, &\text{ else}
                         \end{cases}~.$$
Similarly, the next block sends 
$$\e_{j-1} \mapsto \begin{cases} \e_{j-1}, &\text{ if } j-1>s_{i-2} \\
                                                -\e_{j-2}, &\text{ else}
                         \end{cases}~.$$
Letting $\beta_S(i,j):= \card\{1 \le \ell <i ~|~ j+1-\ell \le s_{i-\ell}\}$, we see that $\widehat{\rot}_S$ sends
$\e_j \mapsto (-1)^\beta \e_{j-\beta}$ and $\e_{j+1} \mapsto \e_{s_i}.$ Therefore, the composite
$\frac{\partial{\trot_S}}{\partial \theta_{(i,j)}}(\Theta_{S<T}^0)$ is given by
$\e_i \mapsto (-1)^{\beta+1} \e_{j-\beta}$ and $\e_{j+k-i+1} \mapsto (-1)^{k-i+1} \e_{s_i}$. The $(i,j)$
column of
$$D_{\Theta^0_{S<T}} \trot_S = \left[ \left(\frac{\partial{\trot_S}}{\partial \theta_{(i,j)}}\right)_{xy} \right]$$
has a single nonzero entry of $(-1)^{\beta+1}$ in the $(j-\beta, i) \in \mc{M}_S$ row. To compute the
determinant of $D_{\Theta^0_{S<T}} \trot_S$, we first observe that for $(i,j)<(i,j-1) \in Z_S$, there is an inequality, 
$j-1-\beta_S(i,j-1)<j-\beta_S(i,j)$. Therefore, ignoring signs, $D_{\Theta^0_{S<T}} \trot_S$, is a block sum of 
antidiagonal matrices. There will be a block of size $s_i-i$ for each $1 \le i \le k$, so ignoring signs, the determinant
of $D_{\Theta^0_{S<T}}\trot_S$ is
$$(-1)^{\sum\limits_{i=1}^k\sum\limits_{r=1}^{s_i-i-1} r}~.$$
The number of negative entries of $D_{\Theta^0_{S<T}}$ is given by
$$\sum\limits_{(i,j)\in Z_S} (\beta_S(i,j)+1) = d(S) + \sum\limits_{(i,j)\in Z_S} \beta_S(i,j)~.$$
All told,
$$\det(D_{\Theta^0_{S<T}}\trot_S) = (-1)^{\sum\limits_{i=1}^k\sum\limits_{r=1}^{s_i-i-1} r}
    (-1)^{d(S)+ \sum\limits_{(i,j)\in Z_S} \beta_S(i,j)}~.$$
The only difference in computing $\det(D_{\Theta^\pi_{S<T}}\trot_S)$ 
is that the matrix in the $(i_{S<T},S_{i_{S<T}})$ slot is no longer the identity matrix, but rather the diagonal 
matrix all of whose diagronal entries are $1$, except the $(i_{S<T},i_{S<T})$ and $(i_{S<T}+1,i_{S<T}+1)$
entries are both $-1$. For $(i,j) \in Z_S$, if $(i,j)<(i_{S<T}, s_{i_{S<T}})$, then $i<i_{S<T}$, and thus 
$\overline{\rot}_S$ still sends $e_i \mapsto e_j$ as in the $\Theta_{S<T}^0$ case. 
Let us now consider the case $(i,j)>(i_{S<T},s_{i_{S<T}})$. If $i=i_{S<T}$, then the effect on
$\widehat{\rot}_S$ is that the $i_{S<T}$ block 
now sends $\e_{j+1} \mapsto -\e_{s_{i_{S<T}}}$ and $\e_{i_{S<T}+1} \mapsto -\e_{i_{S<T}+1}$. Thus,
this block sends $\e_j \mapsto \e_j$ as in the $\Theta^0_{S<T}$ case. The last case to consider is that
$i>i_{S<T}$. The effect on $\widehat{\rot}_S$ is that the $i_{S<T}$ block now sends 
$\e_{i_{S<T}} \mapsto -\e_{i_{S<T}}$ and $\e_{i_{S<T}+1} \mapsto -\e_{i_{S<T}+1}$. This will introduce an
extra factor of $-1$ on the image $\e_j$ precisely if $j+1-(i-i_{S<T}) = s_{i_{S<T}}+1$, 
or more simply if $j-i=s_{i_{S<T}}-i_{S<T}$. Thus, the determinant of $D_{\Theta^\pi_{S<T}}\trot_S$ will have
one extra factor of $-1$ for each pair $(i,j) \in Z_S$ for which $j-i=s_{i_{S<T}}-i_{S<T}$. Note that
$$\card\{(i,j) \in Z_S ~|~ j-i=s_{i_{S<T}}-i_{S<T}\} = k-i_{S<T}~.$$
Therefore, 
$$\det(D_{\Theta^\pi_{S<T}}\trot_S) = (-1)^{\sum\limits_{i=1}^k\sum\limits_{r=1}^{s_i-i-1} r}
    (-1)^{d(S)+\sum\limits_{(i,j)\in Z_S} \beta_S(i,j)}(-1)^{k-i_{S<T}}~.$$
Recall that $\msf{swap}:(0,\pi)^{Z_S} \times \{0\} \times (0,\pi)$ simply moves the last coordinate to the 
$(i_{S<T},s_{i_{S<T}})$ coordinate in $Z_T$. Thus,
the derivative of $\msf{swap}$ is a matrix consisting of $1$'s along the diagonal until the 
$(i_{S<T},s_{i_{S<T}})$ row.
The $(i_{S<T},s_{i_{S<T}})$ row will consist of all zeros, except the final column will be a $1$. There will be $1$'s
along the subdiagonal, and all other entries are 0. Likewise, the derivative of $\msf{Rswap}$ will be the same as
the derivative of $\msf{swap}$, except the last column of the $(i_{S<T},s_{i_{S<T}})$ row will be a $-1$, since
$\msf{Rswap}$ is orientation reversing in that factor. Thus, 
$$\det(D_{\Theta_{S<T}^0}\msf{swap}) = (-1)^{\card\{(i,j) \in Z_T ~|~ (i,j)>(I_{S<T},s_{I_{S<T}})\}}~,$$
and 
$$\det(D_{\Theta_{S<T}^0}\msf{Rswap}) = (-1)^{\card\{(i,j) \in Z_T ~|~ (i,j)>(I_{S<T},s_{I_{S<T}})\}+1}~.$$

All told, we have
\begin{align*}\sigma^T_S([\Theta_{S<T}^0]) &= (-1)^{\sum\limits_{i=1}^k\sum\limits_{r=1}^{s_i-i-1} r + 
    d(S)+\sum\limits_{(i,j)\in Z_S} \beta_S(i,j)+ \sum\limits_{i=1}^k\sum\limits_{r=1}^{t_i-i-1} r +
    d(T)+\sum\limits_{(i,j)\in Z_T} \beta_T(i,j)+ \card\{(i,j) \in Z_T ~|~ (i,j)>(i_{S<T},s_{i_{S<T}})\} } \\
    &= (-1)^{1+\sum\limits_{i=1}^k\sum\limits_{r=1}^{s_i-i-1} r + 
    \sum\limits_{(i,j)\in Z_S} \beta_S(i,j)+ \sum\limits_{i=1}^k\sum\limits_{r=1}^{t_i-i-1} r +
    \sum\limits_{(i,j)\in Z_T} \beta_T(i,j)+ \card\{(i,j) \in Z_T ~|~ (i,j)>(i_{S<T},s_{i_{S<T}})\} }~,
\end{align*}
since $d(S)+d(T) = 2d(S)+1$.
Notice that for $i \neq i_{S<T}$,
$$\sum\limits_{r=1}^{t_i-i-1} r = \sum\limits_{r=1}^{s_i-i-1} r~.$$
Further, $t_{i_{S<T}} = s_{i_{S<T}}+1$.
Thus,
$$(-1)^{\sum\limits_{i=1}^k\sum\limits_{r=1}^{s_i-i-1} r +\sum\limits_{i=1}^k\sum\limits_{r=1}^{t_i-i-1} r}
     = (-1)^{s_{i_{S<T}}-i_{S<T}}~.$$
Next, note that
$$\card\{(i,j) \in Z_T ~|~ (i,j)>(i_{S<T},s_{i_{S<T}})\} = \sum\limits_{i=i_{S<T}}^k (s_i-i)~.$$
To assess when $\beta_S(i,j) \neq \beta_T(i,j)$, we must figure out if there exists $1 \le \ell <i$ for which
$j+1-\ell \le t_{i-\ell}$, yet $j+1-\ell > s_{i-\ell}$. Since $s_i=t_i$ for all $i$ except $i=i_{S<T}$, the only $\ell$
for which this could possibly hold is $\ell=i-i_{S<T}$. Then, we are seeking $(i,j) \in Z_S$ for which
$j-i+i_{S<T}+1 \le t_{i_{S<T}} = s_{i_{S<T}}+1$ and $j-i+i_{S<T}+1 >s_{i_{S<T}}$. Both of these inequalities
hold precisely if $j-i=s_{i_{S<T}}-i_{S<T}$. Therefore,
$$(-1)^{\sum\limits_{(i,j)\in Z_S}\beta_S(i,j) + \sum\limits_{(i,j) \in Z_T}\beta_T(i,j)} = 
    (-1)^{\beta_T(i_{S<T},s_{i_{S<T}})+\card\{(i,j)\in Z_S ~|~ j-i=s_{i_{S<T}}-i_{S<T}\}} =
    (-1)^{k-i_{S<T}}~,$$
since $\beta_T(i_{S<T},s_{i_{S<T}}) = 0$, and $\card\{(i,j)\in Z_S ~|~ j-i=s_{i_{S<T}}-i_{S<T}\}=k-i_{S<T}$.
Hence, the formula for $\sigma^T_S([\Theta_{S<T}^0])$ is proven. The above reductions also yield the stated
formula for $\sigma^T_S([\Theta_{S<T}^\pi])$.
\qed \\

\begin{observation}
    We can further simplify the boundary formula
    \begin{align*}
        \sigma_S^T([\Theta_{S<T}^0])+\sigma_S^T([\Theta_{S<T}^\pi]) 
        &= \sum (-1)^{1+k+i_{S<T}+\sum\limits_{i=i_{S<T}+1}^k s_i-i}+ 
                (-1)^{\sum\limits_{i=i_{S<T}}^k s_i-i} \\
        &= (-1)^{i_{S<T} + \sum_{i=i_{S<T}+1}^k (s_i-i)} \left( (-1)^{(s_{i_{S<T}})} - (-1)^k\right) \\
        &= (-1)^{\sum_{i=i_{S<T}}^k (s_i-i)}\left(1-(-1)^{k-s_{i_{S<T}}}\right) \\
        &= \begin{cases}
                (-1)^{\sum_{i=i_{S<T}}^k (s_i-i)} 2, &\text{ if } S<T \text{ and } k \not\equiv s_{i_{S<T}} \mod 2 \\
                0, &\text{ otherwise}
            \end{cases}~.
    \end{align*}
\end{observation}

\begin{observation}[Cohomology]\label{cohomology}
    Lemma \ref{boundarySchubert} gives a simple easy to compute formula for the boundary maps in 
    $C_*(\gr_k(n);\bb{Z})$, the integral chain complex of the Schubert CW structure on $\gr_k(n)$. Note that this 
    gives us a description of $C^*(\gr_k(n);\bb{Z})$, the integral chain complex. The $i$-th chain group is
    generated by the same generators of $i$-th chain group, and now the differential
    $C^i(\gr_k(n);\bb{Z}) \to C^{i+1}(\gr_k(n);\bb{Z})$ is simply the transpose of the differential
    $C_{i+1}(\gr_k(n);\bb{Z}) \to C_i(\gr_k(n);\bb{Z})$. In other words, given $S<T \in \nchoosek{n}{k}$ with
    $d(T)=d(S)+1$, the coefficient $\sigma^S_T = \sigma^T_S$. 
    For each $0\leq r \leq k$, consider the map
    $$d_r : \nchoosek{n}{k} \xrightarrow{\{s_1<\cdots<s_k\} \mapsto\sum_{r\leq i \leq k} s_i-i} \ZZ_{\geq 0}~.$$
    The differential $\partial^i: C^i(\gr_k(n);\bb{Z}) \to C^{i+1}(\gr_k(n);\bb{Z})$ is given by
        \begin{equation}
            \partial^i : S = \{s_1< \cdots < s_k\} ~\mapsto~ \sum\limits_{r \in \left\{1 \leq r \leq k ~\mid~ s_{r+1}-s_r >1 
                \text{ and } k-s_r \text{ is odd}\right\}} (-1)^{d_{r-1}(S)} 2 \cdot S_r~,
        \end{equation}
        where $S_r:= \{s_1 < \cdots < s_{r-1} < s_r +1 < s_{r+1} <\cdots< s_k \} \in \nchoosek{n}{k}$.
\end{observation}

\begin{observation}[Arbitrary coefficients]\label{coefficients}
    Thus, for the chain complex with coefficients in a commutative ring $R$, we see that 
    $C^*(\gr_k(n); R)$ has the same underlying chain groups, and differential specified by 
    $$\sigma^S_T([\Theta_{S<T}^0])+\sigma^S_T([\Theta_{S<T}^\pi])
        = \begin{cases}
            (-1)^{\sum_{i=i_{S<T}}^k (s_i-i)} 2, &\text{ if } S<T \text{ and } k \not\equiv s_{i_{S<T}} \mod 2 \\
            0, &\text{ otherwise}
        \end{cases}~.$$
    Thus, if $0=2$ in $R$, the differentials are all $0$, which in particular, recovers the case of $\bb{Z}/2\bb{Z}$
    coefficients, e.g. \cite{characteristic}.
    These computations should also specialize to agree with those of \cite{borel} and \cite{takeuchi}, though we do
    not verify this.
\end{observation}

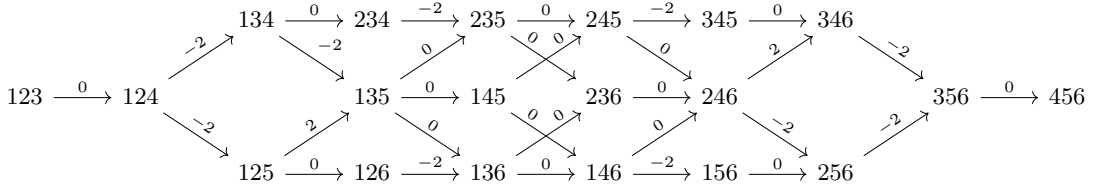
\begin{figure}
\begin{small}\begin{center}
\begin{tikzcd}
&&134\arrow[r, sloped,"0"]\arrow[dr,"-2"]&234\arrow[r,"-2"]&235\arrow[r,"0"]\arrow[dr,near start,"0",sloped]&245\arrow[r,"-2"]\arrow[dr,"0",sloped]&345\arrow[r,sloped,"0"]&346\arrow[dr,"-2",sloped]&&\\
123\arrow[r,"0"]&124\arrow[ur,"-2",sloped]\arrow[dr,"-2",sloped]&&135\arrow[ur, sloped,"0"]\arrow[dr,"0",sloped]\arrow[r,"0"]&145\arrow[ur,near end,"0",sloped]\arrow[dr,near start,"0",sloped]&236\arrow[r,"0"]&246\arrow[ur,"2",sloped]\arrow[dr,"-2",sloped]&&356\arrow[r,"0"]&456\\
&&125\arrow[r,"0"]\arrow[ur, sloped,"2"]&126\arrow[r,"-2"]&136\arrow[ur,near end,sloped,"0"]\arrow[r,"0"]&146\arrow[ur,sloped,"0"]\arrow[r,"-2"]&156\arrow[r,"0"]&256\arrow[ur,"-2",sloped]&&
\end{tikzcd}
\end{center}\end{small}
\caption{A visual of the chain complex for $\gr_3(6)$. Underlying this figure is the poset $\nchoosek{6}{3}$,
with each relation arrow labeled by the coefficient of the differential.}
\label{posetCoho}
\end{figure}

\section{The cohomology of $\gr_k(n)$}\label{sec:cohomology}
The visualization of $C^\ast(\gr_6(12); \bb{Z})$ from Figures \ref{gr510} and \ref{gr612} suggests the chain complex
can be written as a finite direct sum of cubes. 
We will now prove this remarkable observation in Theorem \ref{cochainDecomp}.
A consequence of this theorem is a closed formula for the $R$-cohomology of $\gr_k(n)$,
which is proven in Corollary \ref{cohomologyFormula}. \\

For $S \in \nchoosek{n}{k}$, define
$$\msf{In}(S):= \{i \in \ul{k} ~|~ s_i \equiv k~(\text{mod } 2) \text{ and } s_{i-1} < s_i-1 \} \subset \ul{k}~,$$
and
$$\msf{Out}(S) :=\{i \in \ul{k} ~|~ s_i \not\equiv k ~(\text{mod } 2) \text{ and } s_i+1 < s_{i+1} \} \subset \ul{k}~,$$
where we set $s_0:=0$, and $s_{k+1}:=n+1$, for notational purposes. Further, define
$$\nchoosek{n}{k}_\msf{Out} := \{S \in \nchoosek{n}{k} ~|~ \msf{Out}(S) = \emptyset\}~, 
    \qquad \nchoosek{n}{k}_\msf{In} := \{S \in \nchoosek{n}{k} ~|~ \msf{In}(S) = \emptyset\}~,$$
$$\nchoosek{n}{k}_\msf{Min} := \{S \in \nchoosek{n}{k} ~|~ \msf{Min}(\In(S) \cup \Out(S)) \in \In(S)\}.$$

Consider the subposet $\nchoosek{n}{k}^\adm \subset \nchoosek{n}{k}$ consisting of the same objects as
$\nchoosek{n}{k}$, yet only those relations $S < T$ that factor as a sequence of relations
$S=U_0 < U_1 < \cdots < U_\ell = T$ in which for all $0<r \le \ell$, we have
$d(U_r)-d(U_{r-1}) = 1$ and $i_{U_{r-1}<U_r} \not\equiv k ~\msf{mod}~ 2$. 
This subposet is generated by the relations $S < T$ for which $d(T)-d(S) = 1$ and 
$\sigma^S_T([\Theta_{S<T}^0])+\sigma^S_T([\Theta_{S<T}^\pi]) \neq 0$. 
\begin{observation}\label{admissible}
    Let $S<T$ be a relation in $\nchoosek{n}{k}$. This relation belongs to $\nchoosek{n}{k}^\adm$ if and
    only if $s_i=t_i-1$ for all $i \in \Out(S)\cap \In(T)$ and $s_i=t_i$ for all $i \not\in \Out(S) \cap \In(T)$. 
\end{observation}

\begin{theorem}\label{cochainDecomp}
    Let $R$ be a commutative ring. There is an isomorphism of chain complexes
    \begin{align*}
        C_*(\gr_k(n);\ZZ) &\cong \bigoplus_{S \in \nchoosek{n}{k}_\Out}
            \cone(\ZZ \xrightarrow{2} \ZZ)^{\otimes \In(S)}[d(S)-\card(\In(S))] \\
        &\cong \left(\bigoplus_{S \in \nchoosek{n}{k}_\Out \setminus \nchoosek{n}{k}_\In} \ZZ[d(S)]\right) \oplus
            \left(\bigoplus_{S \in \nchoosek{n}{k}_\msf{Min}} \cone(\ZZ \xrightarrow{2}\ZZ)[d(S)-1]\right)
    \end{align*}
\end{theorem}
{\it Proof.}
The first isomorphism is immediate from Lemmas \ref{Sub} and \ref{cone}. 
The first summand of the second isomorphism comes from the
$S \in \nchoosek{n}{k}$ such that $\In(S) = \emptyset$. The second factor follows from \ref{directSum}.
\qed \\

See Figure \ref{gr510} for a depiction of $C^*(\gr_5(10);R)$ and Figure \ref{gr510WL} for a depiction of the 
isomorphic complex
$$\left(\bigoplus_{S \in \nchoosek{10}{5}_\Out \setminus \nchoosek{n}{k}_\In} \ZZ[d(S)]\right) \oplus
        \left(\bigoplus_{S \in \nchoosek{10}{5}_\msf{Min}} \cone(\ZZ \xrightarrow{2}\ZZ)[d(S)-1]\right)~.$$

We now prove the lemmas that yield the proof of Theorem \ref{cochainDecomp}.
\begin{lemma}\label{Sub}
    There is an isomorphism of chain complexes 
    $$C_*(\gr_k(n);\ZZ) \cong \bigoplus_{S \in \nchoosek{n}{k}_\msf{Out}} C_*(\msf{Sub}(\msf{In}(S))^\op;\ZZ)~.$$
\end{lemma}
{\it Proof.}
We can regard 
$\nchoosek{n}{k}^\adm$ as a weighted level graph with weights given by the coefficients
$\sigma^S_T([\Theta_{S<T}^0])+\sigma^S_T([\Theta_{S<T}^\pi])$, and the level given by the dimension of
the corresponding Schubert cells. Recall there is a bijection between weighted level graphs whose adjacency
matrices square to the zero matrix and chain complexes. Since $C_*(\gr_k(n))$ is a chain complex, and 
$\nchoosek{n}{k}^\adm$ precisely selects out the relations $S < T$ with 
$\sigma^S_T([\Theta_{S<T}^0])+\sigma^S_T([\Theta_{S<T}^\pi]) \neq 0$, the adjacency matrix of 
$(\nchoosek{n}{k}^\adm)^\op$ squares to the zero matrix. Let us denote the chain complex associated to 
$(\nchoosek{n}{k}^\adm)^\op$ by $ C_*((\nchoosek{n}{k}^\adm)^\op)$.  Thus, by construction, 
$C_*(\gr_k(n)) \cong C_*((\nchoosek{n}{k}^\adm)^\op)$.
Next, we claim the canonical functor between posets
$$\coprod_{S \in \nchoosek{n}{k}_\Out} \nchoosek{n}{k}^\adm_{\le S} \to \nchoosek{n}{k}^\adm$$
is an isomorphism.  To prove this, note that for each $S \in \nchoosek{n}{k}_\Out$, the canonical inclusion
$\nchoosek{n}{k}^\adm_{\le S} \hookrightarrow \nchoosek{n}{k}^\adm$ is fully faithful. Thus, to prove the
claim, it is enough to show that for each $T \in \nchoosek{n}{k}^\adm$, there exists a unique
$S \in \nchoosek{n}{k}_\Out$ such that $T \in \nchoosek{n}{k}^\adm_{\le S}$. So, for
$T \in \nchoosek{n}{k}^\adm$, define a set $S^T \in \nchoosek{n}{k}^\adm$ by defining 
$S^T_i = T_i$ in $i \not\in \Out(T)$, and $S^T_i = T_i+1$ otherwise. By construction, 
$S \in \nchoosek{n}{k}_\Out$. Observation \ref{admissible} implies that $T \in \nchoosek{n}{k}^\adm_{\le S^T}$,
and further that $S^T$ is the unique element of $\nchoosek{n}{k}_\Out$ for which this is true. 
Now, observe that for $S \in \nchoosek{n}{k}_\Out$,  there is an isomorphism of posets
$\nchoosek{n}{k}_{\le S}^\adm \xrightarrow{\cong} \Sub(\In(S))$ given by sending $T \mapsto \Out(T)$. 
Therefore, 
\begin{align*}
    C_*(\gr_k(n)) &\cong C_*\left((\nchoosek{n}{k}^\adm)^\op\right) \\
    &\cong C_*\left(\coprod_{S \in \nchoosek{n}{k}_\Out} (\nchoosek{n}{k}^\adm_{\le S})^\op\right) \\
    &\cong C_*\left(\coprod_{S \in \nchoosek{n}{k}_\Out} \Sub(\In(S))^\op\right) \\
    &\cong \bigoplus_{S \in \nchoosek{n}{k}_\Out} C_*\left( \Sub(\In(S))^\op\right)~.
\end{align*}
\qed \\

\begin{lemma}\label{cone}
    There is aa isomorphism of chain complexes
    $$C_*(\Sub(\In(S)^\op) \cong \cone\left(\ZZ \xrightarrow{2} \ZZ\right)^{\otimes \In(S)}~.$$
\end{lemma}
{\it Proof.}
In terms of their associated weighted level graphs, note that both underlying digraphs are isomorphic to the poset
$\Sub(\{1< \cdots < \card(\In(S))\})$. Further, the weights in both weighted level graphs are all $\pm 2$. 
We will now prove that for any cardinality $r \in \bb{Z}_{\ge 0}$, any two choices of weights drawn from the
set $\{\pm2\}$ on the
digraph $\Sub(\ul{r})$, such that their adjacency matrices square to 0, are isomorphic. In fact, proceeding
by induction on $r$, we only need to prove the case $r=2$, as the cases $r=0,1$ are clear. 
Namely, we will show that given any solid diagrams with weights $u_*^*,v_*^* \in \{\pm2\}$, 
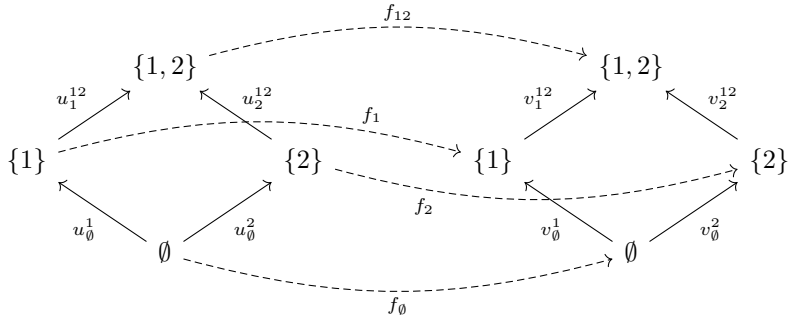
$$\begin{tikzcd}
    & \{1,2\} \arrow[rrrr,dashed,bend left=15,"f_{12}"] & & & & \{1,2\} & \\
    \{1\} \arrow[ur,"u_1^{12}"] \arrow[rrrr,dashed,bend left=15,near end,"f_1"] & & \{2\} \arrow[ul,"u_2^{12}"'] 
        \arrow[rrrr,dashed,bend right=15,near start,"f_2"'] & & \{1\} \arrow[ur,"v_1^{12}"] & & \{2\} 
        \arrow[ul,"v_2^{12}"'] \\
    & \emptyset \arrow[ul,"u_\emptyset^1"] \arrow[ur,"u_\emptyset^2"'] 
        \arrow[rrrr,dashed,bend right=15,"f_\emptyset"'] 
        & & & & \emptyset \arrow[ul,"v_\emptyset^1"] \arrow[ur,"v_\emptyset^2"']& \\
\end{tikzcd}$$
there exists filler arrows, $f_\emptyset,f_1,f_2,f_{12} \in \{\pm1\}$. Further, $f_1,f_2$, and $f_{12}$ are uniquely
determined by $f_\emptyset$. The commutivity conditions give the follwing four equations:\\
\begin{subequations}
    \begin{tabularx}{\textwidth}{Xp{1cm}X}
        \begin{equation}
            \label{c1}
            f_1u_\emptyset^1 = v_\emptyset^1f_\emptyset
        \end{equation}
        & &
        \begin{equation}
            \label{c2}
            f_2u_\emptyset^2 = v_\emptyset^2f_\emptyset
        \end{equation}
    \end{tabularx}
\end{subequations}
\begin{subequations}
    \begin{tabularx}{\textwidth}{Xp{1cm}X}
        \begin{equation}
            \label{c3}
            f_{12}u_1^{12} = v_1^{12}f_1
        \end{equation}
        & &
        \begin{equation}
            \label{c4}
            f_{12}u_2^{12} = v_2^{12}f_2~.
        \end{equation}
    \end{tabularx}
\end{subequations}
The fact that each solid diagram defines a chain complex results in the following two equations: \\
\begin{subequations}
    \begin{tabularx}{\textwidth}{Xp{1cm}X}
        \begin{equation}
            \label{c5}
            u_1^{12}u_\emptyset^1=u_2^{12}u_\emptyset^2
        \end{equation}
        & &
        \begin{equation}
            \label{c6}
            v_1^{12}v_\emptyset^1=v_2^{12}v_\emptyset^2~.
        \end{equation}
    \end{tabularx}
\end{subequations}
Equations \ref{c1} and \ref{c2} uniquely determine $f_1$ and $f_2$ in terms of $f_\emptyset$. 
Multplying Equation \ref{c1} on the left by $v_1^{12}$ and substituting in Equations \ref{c3} and \ref{c6}
yields $f_{12}u_1^{12}u_\emptyset^1=v_2^{12}v_\emptyset^2f_\emptyset$, which uniquely determines
$f_{12}$ in terms of $f_\emptyset$. It remains to show that this value of $f_{12}$ is compatible with the
value determined by Equation \ref{c4}. Multiplying Equation \ref{c2} on the left by $v_2^{12}$ and 
substituting in Equations \ref{c4} and \ref{c5} also yields the equation
$f_{12}u_1^{12}u_\emptyset^1=v_2^{12}v_\emptyset^2f_\emptyset$. Thus, there is a consistent value
of $f_{12}$ uniquely determined by $f_\emptyset$. 
\qed \\

\begin{lemma}\label{directSum}
    For $r >0$, there is an isomorphism of chain complexes
    $$\cone(\ZZ\xrightarrow{2}\ZZ)^{\otimes r} \cong \bigoplus_{0 \le a \le r-1} 
        \left(\cone(\ZZ\xrightarrow{2}\ZZ)[a]\right)^{\oplus {r-1 \choose a}}~.$$
\end{lemma}
{\it Proof.}
We will prove this by induction on $r$. As $r=1$ is clear, consider the case $r=2$. The following diagram
exhibits such an isomorphism:
$$\begin{tikzcd}[ampersand replacement=\&]
    \ZZ \arrow[r,"\begin{bmatrix} 2 \\ 2\end{bmatrix}"] \arrow[dd,"1"] \& \ZZ \oplus \ZZ 
        \arrow[dd,"{\begin{bmatrix} 1 & 0 \\ -1 & 1 \end{bmatrix}}"]
        \arrow[r,"{\begin{bmatrix} 2 & -2\end{bmatrix}}"]  
        \& \ZZ \arrow[dd,"-1"] \\
    \& \& \\
    \ZZ 
        \arrow[r,"\begin{bmatrix} 2 \\ 0 \end{bmatrix}"'] 
        \& \ZZ \oplus \ZZ
        \arrow[r,"{\begin{bmatrix} 0 & 2 \end{bmatrix}}"'] 
        \& \ZZ
\end{tikzcd}.$$
Assume $r>1$, then
\begin{align*}
    \cone(\ZZ \xrightarrow{2} \ZZ)^{\otimes r} &= \cone(\ZZ \xrightarrow{2} \ZZ) \otimes 
        \cone(\ZZ \xrightarrow{2} \ZZ)^{\otimes r-1} \\
    &\cong \cone(\ZZ \xrightarrow{2} \ZZ) \otimes \left( \bigoplus_{0\le a\le r-2} 
        \cone(\ZZ \xrightarrow{2} \ZZ)[a]^{\oplus {r-2 \choose a}}\right) \\
    &\cong \bigoplus_{0\le a \le r-2} \left(\cone(\ZZ \xrightarrow{2} \ZZ) \otimes 
        \cone(\ZZ \xrightarrow{2} \ZZ)\right)[a]^{\oplus{r-2 \choose a}} \\
    &\cong \bigoplus_{0 \le a \le r-2} \left( \cone(\ZZ \xrightarrow{2} \ZZ) \oplus \cone(\ZZ \xrightarrow{2} \ZZ)\right)
        [a]^{\oplus{r-2 \choose a}} \\
    &\cong \left( \bigoplus_{0 \le a \le r-2} \cone(\ZZ \xrightarrow{2} \ZZ)[a]^{\oplus{r-2 \choose a}}\right) \oplus
        \left(\bigoplus_{0 \le a \le r-2} \cone(\ZZ \xrightarrow{2} \ZZ)[a+1]^{\oplus{r-2\choose a}}\right) \\
    &\cong \cone(\ZZ \xrightarrow{2} \ZZ) \oplus \left( \bigoplus_{0 < a \le r-2} 
        \left(\cone(\ZZ \xrightarrow{2} \ZZ)^{\oplus{r-2\choose a}} \oplus
        \cone(\ZZ \xrightarrow{2} \ZZ)^{\oplus{r-2\choose a-1}}\right)[a]\right) \oplus 
        \cone(\ZZ \xrightarrow{2} \ZZ)[r-1] \\
    &\cong \bigoplus_{0 \le a \le r-1} \cone(\ZZ \xrightarrow{2} \ZZ)[a]^{\oplus{r-1 \choose a}}~.
\end{align*}
\qed \\

We now prove several corollaries of Theorem \ref{cochainDecomp}, that, in particular, 
provide a closed formula for the $R$-cohomology of $\gr_k(n)$.

\begin{corollary}\label{cochains}
    Let $R$ be a commutative ring. There is an isomorphism of chain complexes
    \begin{align*}
        C^*(\gr_k(n);R) \cong &\left(\bigoplus_{S \in \nchoosek{n}{k}_\Out \cap \nchoosek{n}{k}_\In} R[-d(s)]\right) \\
        &\oplus
        \left(\bigoplus_{S \in \nchoosek{n}{k}_\Out \setminus \nchoosek{n}{k}_\In} 
        \bigoplus_{0 \le a<\card(\In(S))} \cone(R\xrightarrow{2}R)^{\oplus{\card(\In(S))-1 \choose a}}[-d(S)+a]\right)~.
    \end{align*}
\end{corollary}
{\it Proof.}
This follows from the prior lemmas, and observing that 
$$\msf{Hom}(\cone(\bb{Z}\xrightarrow{2}\bb{Z}),R) \cong \cone(R\xrightarrow{2}R)[-1]~.$$
\qed \\

\begin{corollary}\label{cohomologyFormula}
    There is an isomorphism of graded $R$-modules
    $$H^*(\gr_k(n);R) \cong \bigoplus_{S \in \nchoosek{n}{k}} V_S[-d(S)]~,$$
    where 
    $$V_S := \begin{cases}
                            R, &\text{ if } \msf{In}(S)  = \emptyset = \msf{Out}(S) \\
                            \ker(R\xrightarrow{2}R), &\text{ if } \msf{Min}(\In(S) \cup \Out(S)) \in \Out(S) \\
                            \coker(R\xrightarrow{2}R), &\text{ if } \msf{Min}(\In(S) \cup \Out(S)) \in \In(S)
                    \end{cases}.$$
\end{corollary}
{\it Proof.} 
Taking the cohomology of the chain complex in Corollary \ref{cochains} yields
\begin{align*}
    H^*(\gr_k(n);R)
    &\cong \left( \bigoplus_{S \in \nchoosek{n}{k}_\Out \cap \nchoosek{n}{k}_\In}
        R[d(S)]\right) \\
    &\oplus\bigoplus_{S \in \nchoosek{n}{k}_\Out \setminus \nchoosek{n}{k}_\In}
        \bigoplus_{0 \le a < \card(\In(S))} H^*\left(0 \to R[a] \xrightarrow{2} R[a-1] \to 0\right)[-d(S)+\card(\In(S))]^{
        \oplus {\card(\In(S))-1 \choose a}} ~.
\end{align*}
Note that
\begin{align*}
    H^*&\left(0 \to R[a] \xrightarrow{2} R[a-1] \to 0\right)[-d(S)+\card(\In(S))] \\
    &\cong
    \ker(R\xrightarrow{2}R)[-d(S)+\card(\In(S))-a] \oplus \coker(R\xrightarrow{2}R)[-d(S)+\card(\In(S))-a-1]~.
\end{align*}
Thus, the number of summands in 
$$\bigoplus_{0 \le a < \card(\In(S))} H^*\left(0 \to R[a] \xrightarrow{2} R[a+1] \to 0\right)[d(S)-\card(\In(S))]^{
        \oplus {\card(\In(S))-1 \choose a}}$$
is equal to $2^{\card(\In(S))}$. Hence there is one summand for each set in
$\nchoosek{n}{k}^\adm_{\le S}$. Let us choose the following convention: Assign to each
$T \in \nchoosek{n}{k}^\adm_{\le S}$, the $R$-module
$$T \mapsto V_T :=\begin{cases} \ker(R\xrightarrow{2}R), &\text{ if } \msf{Min}(\In(T) \cup \Out(T)) \in \Out(T) \\
                                                        \coker(R\xrightarrow{2}R), &\text{ if } \msf{Min}(\In(T) \cup \Out(T)) \in \In(T)
                                \end{cases}.$$
This choice amounts to the following: Consider the weighted level graph assoicatied to the complex in the description of
$C^*(\gr_k(n);R)$ from Corollary \ref{cochains}. See Figure \ref{gr510WL} for the example of $\gr_5(10)$.
We assign $\ker(R\xrightarrow{2}R)$ to the start node of each edge, and $\coker(R\xrightarrow{2}R)$ to the end
node of each edge.
Recall from the proof of Proposition \ref{Sub} that for each $T$, there is a unique $S \in \nchoosek{n}{k}$ such that
$T \in \nchoosek{n}{k}^\adm_{\le S}$. Therefore, 
\begin{align*}
    \bigoplus_{S \in \nchoosek{n}{k}_\Out \setminus \nchoosek{n}{k}_\In}
        \bigoplus_{0 \le a < \card(\In(S))} &H^*\left(0 \to R[a] \xrightarrow{2} R[a-1] \to 0\right)[-d(S)+\card(\In(S))]^{
        \oplus {\card(\In(S))-1 \choose a}} \\
    &\cong \bigoplus_{S \in \nchoosek{n}{k}_\Out \setminus \nchoosek{n}{k}_\In}
        \bigoplus_{T \in \nchoosek{n}{k}^\adm_{\le  S}} V_T[-d(T)] \\
    &\cong \bigoplus_{S \in \nchoosek{n}{k} \setminus\left(\nchoosek{n}{k}_\Out \setminus 
        \nchoosek{n}{k}_\In\right)} V_S[-d(S)]~.
\end{align*}
Extending the assignment $S \mapsto V_S$ by defining $V_S:= R$ if $\In(S)=\emptyset=\Out(S)$, yields the
desired result:
$H^*(\gr_k(n);R) \cong \bigoplus_{S \in \nchoosek{n}{k}} V_S[-d(S)]$.
\qed \\

\section{Computations}\label{sec:computations}
\begin{figure}[H]
\centering
\scalebox{.6}{\input{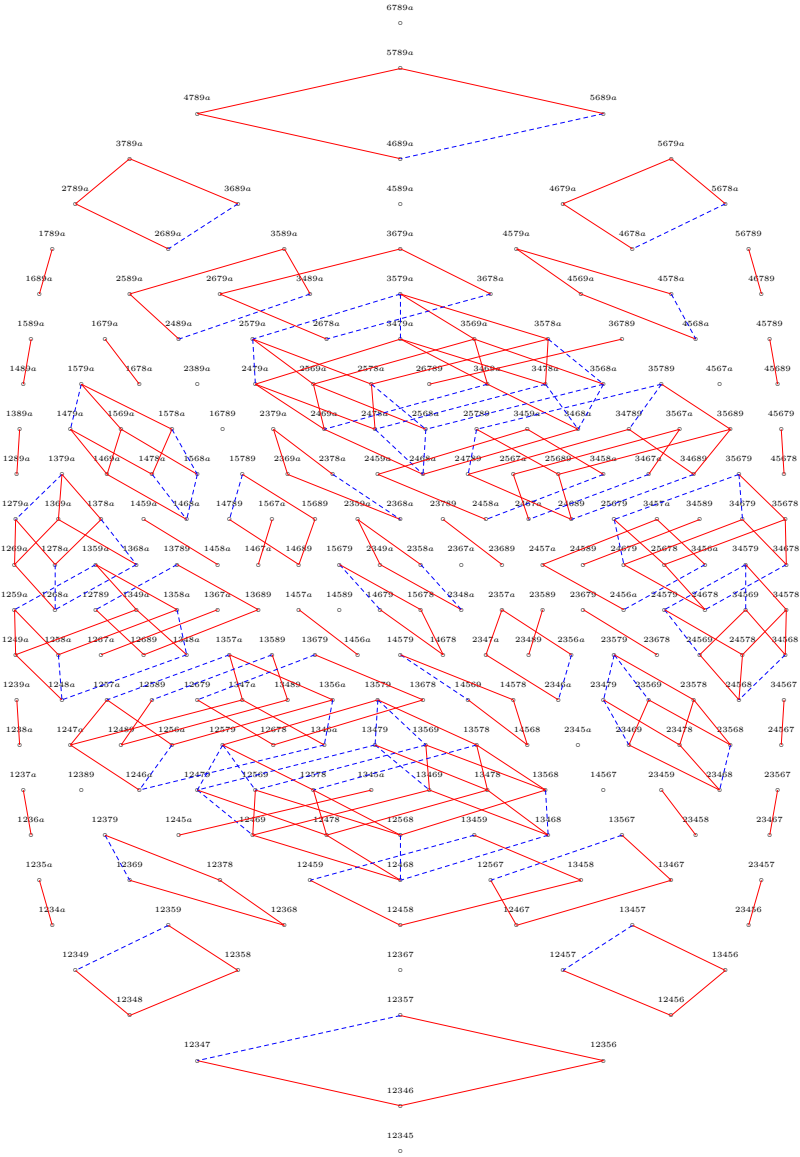}}
\caption{ The chain complex for $\gr_5(10)$. The dotted blue lines indicate the differential is $+2$, the solid
red lines indicate the coefficient of the differential is $-2$. The sets are labeled using base 11: $a = 10$.}
\label{gr510}
\end{figure}

\begin{figure}[H]
\centering
\scalebox{.6}{\input{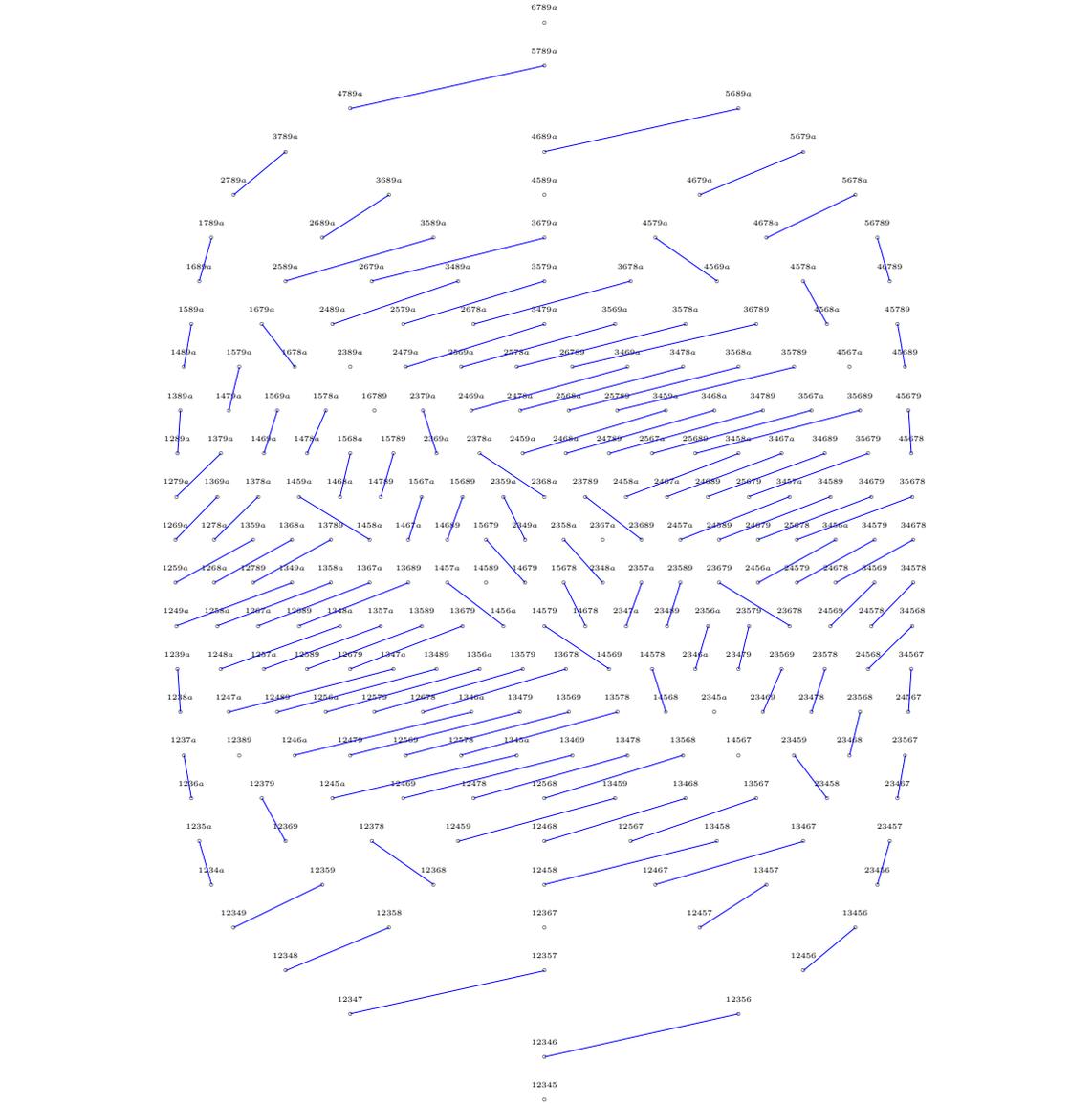}}
\caption{ The chain complex from Theorem \ref{cochainDecomp} that is isomorphic to $C^*(\gr_5(10);R)$. Each blue line corresponds to one of the copies of $\cone(R\xrightarrow{2}R)$. The sets are labeled using base 11: $a = 10$.}
\label{gr510WL}
\end{figure}

\begin{figure}[H]
\centering
\scalebox{.55}{\input{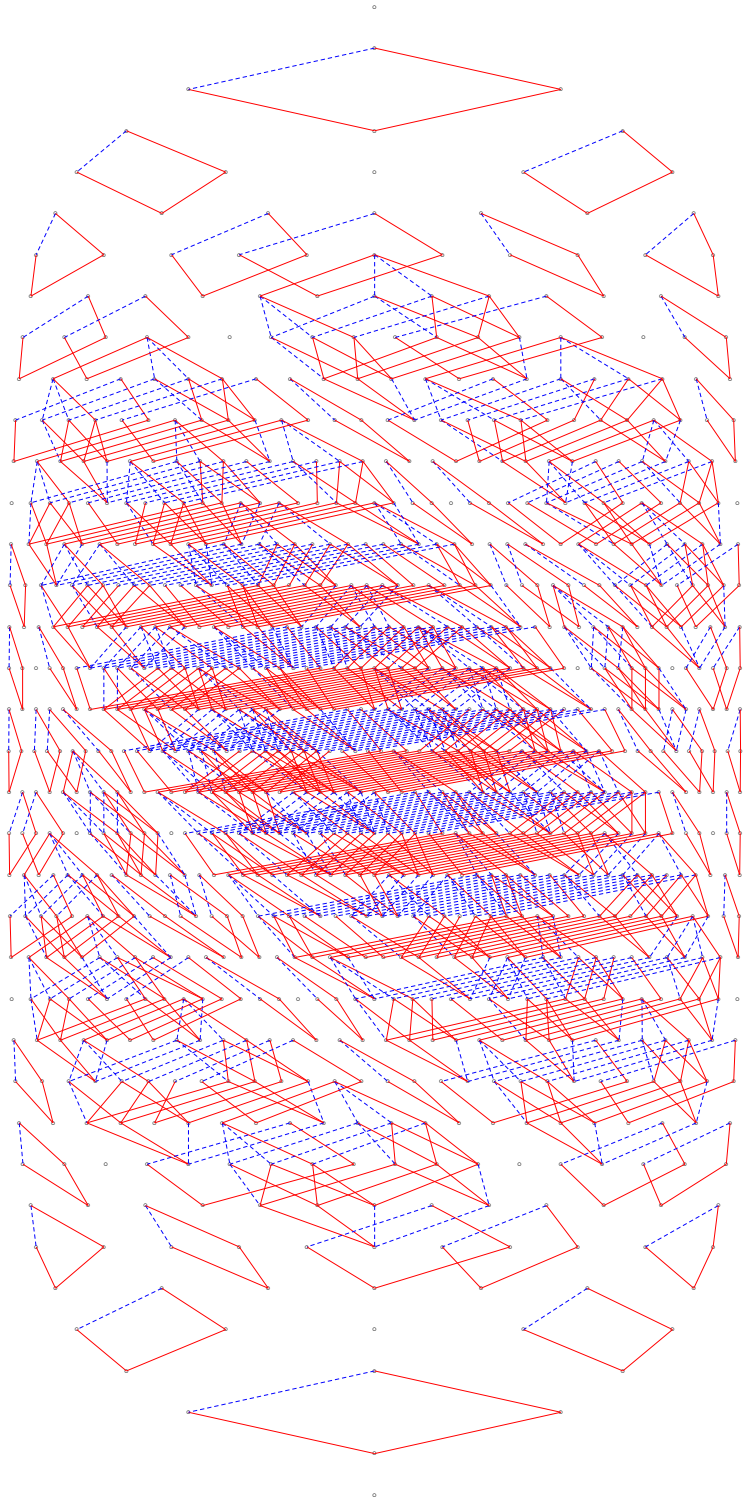}}
\caption{ The chain complex for $\gr_6(12)$. The dotted blue lines indicate the differential is $+2$, the solid
red lines indicate the coefficient of the differential is $-2$. The labels of the nodes were omitted for legibility reasons.}
\label{gr612}
\end{figure}

Finally, we present some computations for various $n$ and $k$.
Figure \ref{posetCoho} provides a graphical depiction of the chain complex of $\gr_3(6)$. The vertical columns
list the generators of $C^i(\gr_3(6))$ in terms of increasing dimension, with the generator of $C^0(\gr_3(6))$ on 
the left, and the generator of $C^6(\gr_3(6))$ on the right. The arrows specify the poset structure of the poset,
$\nchoosek{6}{3}$, which stratifies the space $\gr_3(6)$. Each arrow, $S \to T$, is labeled by the coefficient of its 
differential: $\sigma^S_T([\Theta_{S<T}^0]) + \sigma^S_T([\Theta_{S<T}^\pi])$. 
Figures \ref{gr510} and \ref{gr612} give a graphical depiction of $C^*(\gr_5(10);\bb{Z})$ and 
$C^*(\gr_6(12);\bb{Z})$, respectively. 
Each node in the images represents a set $S \in \nchoosek{n}{k}$. To make the images more legible,
in Figure \ref{gr510}, the nodes have been labeled using base 11 notation, where $a=10$, and we left the nodes
unlabeled in Figure \ref{gr612}.
The blue dotted arrows indicate the coefficient of the differential is $-2$
and the solid red lines indicate the coefficient of the differential is $+2$. 
Figure \ref{gr510WL} is a depiction of the chain complex from Theorem \ref{cochainDecomp} that is isomorphic to $C^*(\gr_5(10);R)$.
Table \ref{gr1224} displays the number of $\Z$ and $\Z/2\Z$ summands in each integral cohomology group
of $\gr_{12}(24)$.

\begin{table}[H]
\centering
\begin{tabular}{|c|c|c|}
\hline
$Gr_{12}(24)$&$\mathbb{Z}$&$\mathbb{Z}/2\mathbb{Z}$\\ \hline
$H^{0}$&1&0\\\hline
$H^{1}$&0&0\\\hline
$H^{2}$&0&1\\\hline
$H^{3}$&0&1\\\hline
$H^{4}$&1&2\\\hline
$H^{5}$&0&2\\\hline
$H^{6}$&0&5\\\hline
$H^{7}$&0&6\\\hline
$H^{8}$&2&9\\\hline
$H^{9}$&0&11\\\hline
$H^{10}$&0&19\\\hline
$H^{11}$&0&23\\\hline
$H^{12}$&3&33\\\hline
$H^{13}$&0&41\\\hline
$H^{14}$&0&58\\\hline
$H^{15}$&0&73\\\hline
$H^{16}$&5&95\\\hline
$H^{17}$&0&117\\\hline
$H^{18}$&0&156\\\hline
$H^{19}$&0&191\\\hline
$H^{20}$&7&239\\\hline
$H^{21}$&0&291\\\hline
$H^{22}$&0&367\\\hline
$H^{23}$&0&441\\\hline
$H^{24}$&11&536\\\hline
$H^{25}$&0&638\\\hline
$H^{26}$&0&777\\\hline
$H^{27}$&0&916\\\hline
$H^{28}$&13&1087\\\hline
$H^{29}$&0&1268\\\hline
$H^{30}$&0&1503\\\hline
$H^{31}$&0&1740\\\hline
$H^{32}$&18&2017\\\hline
$H^{33}$&0&2316\\\hline
$H^{34}$&0&2680\\\hline
$H^{35}$&0&3048\\\hline
$H^{36}$&22&3470\\\hline
$H^{37}$&0&3917\\\hline
$H^{38}$&0&4440\\\hline
$H^{39}$&0&4974\\\hline
$H^{40}$&28&5562\\\hline
\end{tabular}
\begin{tabular}{|c|c|c|}
\hline
$H^{41}$&0&6179\\\hline
$H^{42}$&0&6886\\\hline
$H^{43}$&0&7594\\\hline
$H^{44}$&32&8356\\\hline
$H^{45}$&0&9151\\\hline
$H^{46}$&0&10029\\\hline
$H^{47}$&0&10901\\\hline
$H^{48}$&39&11821\\\hline
$H^{49}$&0&12760\\\hline
$H^{50}$&0&13776\\\hline
$H^{51}$&0&14770\\\hline
$H^{52}$&42&15791\\\hline
$H^{53}$&0&16814\\\hline
$H^{54}$&0&17898\\\hline
$H^{55}$&0&18936\\\hline
$H^{56}$&48&19967\\\hline
$H^{57}$&0&20990\\\hline
$H^{58}$&0&22040\\\hline
$H^{59}$&0&23015\\\hline
$H^{60}$&51&23959\\\hline
$H^{61}$&0&24864\\\hline
$H^{62}$&0&25764\\\hline
$H^{63}$&0&26573\\\hline
$H^{64}$&55&27307\\\hline
$H^{65}$&0&27984\\\hline
$H^{66}$&0&28635\\\hline
$H^{67}$&0&29166\\\hline
$H^{68}$&55&29596\\\hline
$H^{69}$&0&29963\\\hline
$H^{70}$&0&30272\\\hline
$H^{71}$&0&30456\\\hline
$H^{72}$&58&30525\\\hline
$H^{73}$&0&30525\\\hline
$H^{74}$&0&30456\\\hline
$H^{75}$&0&30272\\\hline
$H^{76}$&55&29963\\\hline
$H^{77}$&0&29596\\\hline
$H^{78}$&0&29166\\\hline
$H^{79}$&0&28635\\\hline
$H^{80}$&55&27984\\\hline
$H^{81}$&0&27307\\\hline
$H^{82}$&0&26573\\ \hline
\end{tabular}
\begin{tabular}{|c|c|c|}
\hline
$H^{83}$&0&25764\\\hline
$H^{84}$&51&24864\\\hline
$H^{85}$&0&23959\\\hline
$H^{86}$&0&23015\\\hline
$H^{87}$&0&22040\\\hline
$H^{88}$&48&20990\\\hline
$H^{89}$&0&19967\\\hline
$H^{90}$&0&18936\\\hline
$H^{91}$&0&17898\\\hline
$H^{92}$&42&16814\\\hline
$H^{93}$&0&15791\\\hline
$H^{94}$&0&14770\\\hline
$H^{95}$&0&13776\\\hline
$H^{96}$&39&12760\\\hline
$H^{97}$&0&11821\\\hline
$H^{98}$&0&10901\\\hline
$H^{99}$&0&10029\\\hline
$H^{100}$&32&9151\\\hline
$H^{101}$&0&8356\\\hline
$H^{102}$&0&7594\\\hline
$H^{103}$&0&6886\\\hline
$H^{104}$&28&6179\\\hline
$H^{105}$&0&5562\\\hline
$H^{106}$&0&4974\\\hline
$H^{107}$&0&4440\\\hline
$H^{108}$&22&3917\\\hline
$H^{109}$&0&3470\\\hline
$H^{110}$&0&3048\\\hline
$H^{111}$&0&2680\\\hline
$H^{112}$&18&2316\\\hline
$H^{113}$&0&2017\\\hline
$H^{114}$&0&1740\\\hline
$H^{115}$&0&1503\\\hline
$H^{116}$&13&1268\\\hline
$H^{117}$&0&1087\\\hline
$H^{118}$&0&916\\\hline
$H^{119}$&0&777\\\hline
$H^{120}$&11&638\\\hline
$H^{121}$&0&536\\\hline
$H^{122}$&0&441\\\hline
$H^{123}$&0&367\\\hline
$H^{124}$&7&291\\ \hline
\end{tabular}
\begin{tabular}{|c|c|c|}
\hline
$H^{125}$&0&239\\\hline
$H^{126}$&0&191\\\hline
$H^{127}$&0&156\\\hline
$H^{128}$&5&117\\\hline
$H^{129}$&0&95\\\hline
$H^{130}$&0&73\\\hline
$H^{131}$&0&58\\\hline
$H^{132}$&3&41\\\hline
$H^{133}$&0&33\\\hline
$H^{134}$&0&23\\\hline
$H^{135}$&0&19\\\hline
$H^{136}$&2&11\\\hline
$H^{137}$&0&9\\\hline
$H^{138}$&0&6\\\hline
$H^{139}$&0&5\\\hline
$H^{140}$&1&2\\\hline
$H^{141}$&0&2\\\hline
$H^{142}$&0&1\\\hline
$H^{143}$&0&1\\\hline
$H^{144}$&1&0\\ \hline
\end{tabular}
\caption{This table displays the integral cohomology of $\gr_{12}(24)$. We indicate the number of $\Z$ and $\Z/2\Z$ summands in each cohomology group.}
\label{gr1224}
\end{table}




\bibliography{grassmannians}
\bibliographystyle{alpha}


\end{document}